\documentclass[a4paper, 12pt]{article}
\usepackage[utf8]{inputenc}
\usepackage[T2A]{fontenc}
\usepackage{amsmath,amssymb,amsthm, amsfonts, mathtools, bm,bbm}
\usepackage[a4paper,hmargin=2.5cm,vmargin=2.5cm]{geometry}
\usepackage{xcolor}
\usepackage[english]{babel}

\theoremstyle{definition}

\newtheorem{thm}{Theorem}
\newtheorem{defn}[thm]{Definition}
\newtheorem{rem}{Remark}
\newtheorem{ex}{Example}
\newtheorem{st}[thm]{Proposition}
\newtheorem{lemma}[thm]{Lemma}
\newtheorem{cor}[thm]{Corollary}
\newtheorem*{pf}{Proof}

\newcommand{\BC}{\mathbb{C}}

\newcommand{\BH}{\mathbb{H}}
\newcommand{\Bh}{\mathbbm{h}}
\newcommand{\BI}{\mathbb{I}}

\newcommand{\BM}{\mathbb{M}}

\newcommand{\BS}{\mathbb{S}}
\newcommand{\Bs}{\mathbbm{s}}
\newcommand{\BT}{\mathbb{T}}

\newcommand{\fg}{\mathfrak{g}}

\def\gl{\mathfrak{gl}}

\begin{document}

\makeatletter
\renewcommand{\theequation}{{\thesection}.{\arabic{equation}}}
\@addtoreset{equation}{section} \makeatother

\title{Quantum $\mathfrak{gl}$-weight system and its average values}
\author{\rule{0pt}{7mm} Mikhail Zaitsev \thanks{mrzaytsev@edu.hse.ru}\\
%{\small\it }\\
%{\small\it }\\
%\rule{0pt}{7mm} \thanks{}\\
{\small\it
National Research University Higher School of Economics}\\
%{\small\it }\\
%{\small \it }\\
%{\small \it
%}\\
%{\small \it }
}

\maketitle

\begin{abstract}
We present a proof of a recent conjecture due to M.~Kazarian,
E.~Krasilnikov, S.~Lando, and M.~Shapiro, 
which describes the average value of the universal ${\gl}$-weight system on permutations. The proof uses a quantum analogue of the ${\gl}$-weight system on Hecke algebras of type $A$, which leads to a one-parameter
deformation of the average value of the universal ${\gl}$-weight system. 
We show that the average value of the quantum weight system 
is a linear combination of one-part Schur functions, with
coefficients being $q$-analogues of Bernoulli polynomials.
\end{abstract}

%{\bf AMS Mathematics Subject Classification, 2010:} 05E05, 20C08

{\bf Keywords:} Lie algebra, Hecke algebra, weight system, q-Bernoulli numbers,
Reflection Equation algebra, quantum weight system

\section{Notation}

\subsection{Objects}

\begin{tabular}{cp{0.8\textwidth}}
	$\fg, \langle\cdot,\cdot\rangle$ & a Lie algebra endowed with a nondegenerate invariant bilinear product\\
	$\gl(N)$ & general linear Lie algebra; consists of all $N\times N$ matrices with the commutator serving as the Lie bracket \\
	$d$ & dimension of Lie algebra; specifically, for $\gl(N)$, $d=N^2$ \\
	$E_i^j$, & matrix units, the standard generators of the Lie algebra $\gl(N)$\\
	$1\le i,j\le N$&\\
	$C_1,\cdots,C_N$ & Casimir elements in $U(\gl(N))$\\
	$w$ & a weight system\\
	$w_\fg$ & the Lie algebra weight system associated to a Lie algebra $\fg$\\
	$m$ & the number of permuted elements\\
	$\BS_m$ & symmetric group; \\
	$\alpha$ & a permutation, $\alpha\in\BS_m$\\
	$\BC[\BS_m],\BC(q)[\BS_m]$ & group algebras of the symmetric group\\
	$\Bs_m$ & symmetrizer, $\Bs_m=\frac1{m!}\sum_{\alpha\in\BS_m}\alpha\in\BC[BS_m]$\\
	$\BH_m$ & $m$~th Hecke algebra\\
	$g_i$ & generators of $\BH_m$\\
	$C_{mk}$ & basic elements in $\BH_m$\\
    $\Bh_m$ & quantum symmetrizer, $\Bh_m\in BH_m$ \\
	$\BM(N)$ & Reflection Equation algebra\\
\end{tabular}\\

\subsection{Mappings}

\begin{tabular}{cp{0.8\textwidth}}
	$w_{\gl(N)}$ & the function on permutations associated to the Lie algebra~$\gl(N)$, $w_{\gl(N)}:\BC[\BS_m]\to ZU\gl(N)\equiv \BC[C_1,\dots,C_N]$\\
	$w_{\gl}$ & the  universal $\gl$-weight system, $w_{\gl}:\BC[\BS_m]\to \BC[N;C_1,C_2,\dots]$\\
	$\chi_{m,N}$ & the characteristic mapping, $\chi_{m,N}:\BH_m\to \BM(N)$\\
	$\chi_{m}$ & the universal characteristic mapping, $\chi_{m}:\BH_m\to \BC(q)[p_1,p_2,\dots]$\\
	$\omega_{m,N}$ & the quantum $\gl(N)$-weight system, $\omega_{m,N}:\BH_m\to \BC(q)[C_1,\dots,C_N]$\\
$\omega_{m}$ & the  universal quantum $\gl$-weight system, $\omega_m:\BH_m\to \BC(q)[q^{-N};C_1,C_2,\dots]$\\
$\rho_m$ & the $R$-matrix representation, $\rho_m:\BH_m\to Aut(V^{\otimes m})$\\
$\BT r_{m}$ & the quantum trace on the Hecke algebra, $\BT r_{m}: \BH_m \to \BH_{m-1}$\\	
	\end{tabular}\\

\subsection{Polynomials}

\begin{tabular}{cp{0.8\textwidth}}
	$\nu$ & $\nu=N-1$\\
$B_l(\xi)$ & Bernoulli polynomial\\

 $\beta_{l}^{(h,k)}(\xi)$
 & quantum $l$~th  $q$-Bernoulli polynomial of order $k$\\
 $B_l^{(k)}(\xi)$ &  Bernoulli polynomial of order $k$, $B_l^{(k)}(\xi)=\lim\limits_{q\rightarrow 1}\beta_l^{(h,k)}(\xi)$\\
 	$p_1,\cdots,p_N$ & quantum power sum polynomials;\\
  $W_m(N)$ & value of the ${\frak{gl}}$-weight system on the symmetrizer $\Bs_m$\\
 $\Omega_m(N)$ & value of the quantum ${\frak{gl}}$-weight system on the $q$-symmetrizer $\Bh_m$\\

\end{tabular}\\

\subsection{Matrices}

\begin{tabular}{cp{0.8\textwidth}}
	$E_i^j$ & matrix unit\\
	$E=||E_i^j||$ & matrix of basic elements in $\gl(N)$\\
	$m_i^j$ & generators of $\BM(N)$\\
	$M=||m_i^j||$ & matrix of basic elements in $\BM(N)$\\	
	$R$ & Drinfeld--Jimbo $R$ matrix\\
	$R_{i,i+1}$ & Drinfeld--Jimbo $R$ matrix acting on the tensor product of
	the $i$~th and $i+1$~st factors\\
	$\BI$ & $N\times N$ identity matrix\\
	$P_\alpha$ & the permutation matrix associated to a permutation $\alpha\in\BS_m$, $P_{\alpha} \in End( \mathbb{C}^N)^{\otimes m}$
\end{tabular}\\

\section{Introduction}
The article is devoted to a proof of a conjecture (Theorem \ref{MainThm} below) proposed in \cite{KLUI}.
The conjecture is an explicit expression of 
the average value of the ${\gl}$-weight system 
on permutations as a linear combination of one-part Schur functions. 
Weight systems are functions on chord diagrams, which encode finite
type knot invariants. The extension of the weight systems associated
to the series~$A$ Lie algebras to permutations was suggested 
as a tool to make their computation effective, by means of
a recurrence relation for them~\cite{KL,ZY}. As an important byproduct,
the recurrence allows for unifying all the $\gl(N)$-weight systems,
for $N=1,2,3,\dots$, into a universal $\gl$-weight system on permutations.
The latter takes values in the ring of polynomials in infinitely
many variables $\BC[N][C_1,C_2,\dots]$. 

Our proof of the conjecture in question
makes use of the quantum analogue of the ${\frak{gl}}$-weight system on Hecke algebras of type $A$. It was defined in \cite{GS} using the Reflection Equation algebra. The quantum ${\frak{gl}}$-weight system inherits many properties of the classical ${\frak{gl}}$-weight system. In the proofs of the properties we make substantial use of the Reflection Equation algebra, but all the concepts
and notions we use can be introduced independently of the latter.

The group algebra $\BC[\BS_m]$ of the symmetric group $\BS_n$ is isomorphic
to the Hecke algebra $\BH_m$, for each $m=1,2,\dots$,
but there is no canonical isomorphism between the two.
In particular, we know no natural way to define the quantum
universal $\gl$-weight system on permutations. 
However, the \emph{centers} of the two algebras can be identified canonically,
which makes it possible to define the values of the quantum $\gl$-weight
system on central elements and, in particular, on the central idempotent
entering the conjecture. 

The main result of the present paper is an explicit formula for the mean value of the quantum ${\frak{gl}}$-weight system $\Omega_m$ in terms of one-part Schur functions $S_k$. Set $\nu=N-1$.
Our main result asserts that the value of the averaged quantum ${\frak{gl}}$-weight system is given by the formula 
 \begin{equation} \label{MainMain}
	\Omega_{m}(N) =\frac{m!}{[m]_q!}\frac{[m+\nu]_q!}{(m+\nu)!} \sum_{l=0}^{m} 
	\frac{(m+\nu)_lq^{-2l}}{l!} \beta_{l}^{(m-\nu+1,\nu)}\left(\frac{\nu}{2}\right)  S_{m-l},
\end{equation}
    where, for a positive integer~$b$, we use the standard
    notation $(a)_b$, which denotes the
    falling factorial $(a)_b=a(a-1)(a-2)\dots(a-b+1)$, 
       and
    $\beta_{l}^{(h,k)}(\xi)$ is the \emph{quantum $l$~th $q$-Bernoulli polynomial} of order $k$ (see \cite{C} or \cite{KCR})
$$\beta_{l}^{(h,k)}(\xi) = \frac{1}{(1-q^{-2})^l} \sum\limits_{i=0}^{l}\binom{l}{i} \frac{(i+h)_k}{([i+h]_q)_k }(-q^{-2\xi})^i,$$ 
where we extend the notation for the falling factorial by setting 
$$
([a]_q)_b=[a]_q[a-1]_q\dots[a-b+1]_q.
$$

The quantum $l$~th $q$-Bernoulli polynomial $\beta_l^{(h,k)}(\xi)$ has a limit as $q\to1$. This limit $$B_l^{(k)}(\xi)=\lim\limits_{q\rightarrow 1}\beta_l^{(h,k)}(\xi)$$
does not depend on the value of $h$.  The generating power series
for these limit polynomials is
$$\sum_{l=0}^\infty B_l^{(k)}(\xi)\frac{u^l}{l!}=e^{\xi t}\left(\frac{e^t-1}{t}\right)^{-k}.$$
As a consequence, each summand in the right hand side of 
Eq.~(\ref{MainMain}) has a limit as $q\to 1$ and these limit values are
\begin{equation} \label{AlmostMain}
    \underset{q \rightarrow 1}{\lim}~ \Omega_{m}(N) = \sum_{l=0}^m\frac{(\nu+m)_l}{l!}B_l^{( \nu)}\left(\frac\nu2\right)S_{m-l}.
\end{equation}
We will prove that this limit value coincides with the average value $W_m(N)$
of the ${\frak{gl}}$-weight system on permutations of~$m$ elements. 
Expression (\ref{AlmostMain}) for~$W_m(N)$ was proposed in the article \cite{KLUI}. The coefficients $B_l^{(k)}$ in the expansion are called \emph{Bernoulli polynomials of order $k$}.

The paper is organized as follows. 

In Section~\ref{s-op} we give an outline of the proof of the conjecture,
without giving necessary definitions.

In the first part of Section~\ref{sglN} we 
give an alternative definition of the $\frak{gl}$-weight system 
using matrix notation in the universal enveloping algebra $U({\gl(N)})$. 
In the second part of Section~\ref{sglN}, we introduce a quantum analog of the $\frak{gl}$-weight system defined on the Hecke algebras of type $A$. In the last part of the section, we recall some basic facts about $R$-matrix representations of the Hecke algebra (see the survey paper \cite{OP}) and about the Reflection Equation algebra (see~\cite{GPS}). Using these facts, we construct a model for the quantum $\frak{gl}$-weight system and prove its existence.  

Section~\ref{sqglN} is devoted to the definition of the average value of the quantum $\frak{gl}$-weight system. The existence of a quasi-classical limit $q \rightarrow 1$ is proved.

Section~\ref{sHgl} is entirely dedicated to the proof of formulas (\ref{MainMain}) and (\ref{AlmostMain}). In the first part of Section 4 we describe the main statements of the section and prove formulas (\ref{MainMain}) and (\ref{AlmostMain}) using them. In the second part of the section we calculate the average value of the quantum ${\frak{gl}}$-weight system, and in the last subsection we compute the limit for the coefficients in the expansion.

The author is grateful to M.E. Kazarian and S.K. Lando for their interest in his work, helpful discussions, and substantial contributions to the manuscript.

\section{Outline of the proof}\label{s-op}

Before giving formal definitions we present a brief outline of the
proof of the main result, in order to make the paper more
readable. 

In addition to the $\gl(N)$ weight systems
$$
w_{\gl(N)}:\BC[\BS_m]\to \BC[C_1,C_2,\dots,C_N]
$$
and the universal $\gl$ weight system
$$
w_{\gl}:\BC[\BS_m]\to \BC[N][C_1,C_2,\dots]
$$
we make use of the quantum weight systems
$$
\omega_{m,N}:\BH_m\to \BM(N)
$$
and the universal quantum weight system
$$
\omega_{m}:\BH_m\to \BC[q^{-N}][C_1,C_2,\dots]
$$
defined in~\cite{GS}. The domain $\BH_m$ of the quantum weight systems
is the $m$~th Hecke algebra, and their range is the Reflection Equation
algebra $\BM(N)$.

The universal quantum weight system $\omega_m$ can be defined by
axiomatizing its properties, see Definition~\ref{UnivWeightSystemDef} below.
The sequence $\omega_{m,N}$, $m=1,2,\dots$, of quantum weight systems 
can be defined using an explicit formula in the Reflection Equation algebra, see Definition~\ref{d-qws}
below. They take values in the centers of the Reflection Equation algebras
$\BM(N)$. In turn, this center admits a natural identification with
the ring of polynomials $\BC[C_1,C_2,\dots,C_N]$. The sequence of universal 
quantum weight systems $\omega_m$ possesses the property that it 
specializes into $\omega_{m,N}$, for a given~$N$.

The averaged sum 
$$
\Bs_m=\frac1{m!}\sum_{\alpha\in\BS_m}\alpha\in\BC[\BS_m]
$$
of all permutations in the group~$\BS_m$ is an idempotent in~$\BC[\BS_m]$.
Using the natural inclusion $\BS_{m-1}\subset\BS_m$,
the sequence of these idempotents, for $m=1,2,\dots$, can be defined recursively
by the equation
$$
\Bs_m=\frac1m\Bs_{m-1}(id_m+(m-1)\cdot(m-1,m))\Bs_{m-1},\quad m\ge2,
$$
with the initial condition $\Bs_1=1$.

The sequence of Hecke algebras $\BH_m$,  $m=1,2,\dots$, contains idempotents
$\Bh_m$ defined  in the following similar way:
$$
\Bh_1=1,\qquad\Bh_m=\frac1{[m]_q}\Bh_{m-1}(q^{2-2m}+q^{-1}[m-1]_qg_{m-1}))\Bh_{m-1},
$$
where $g_i$, $i=1,2,\dots$ are standard Artin generators of the Hecke algebras.

Then we manage to prove the following statements:
\begin{itemize}
	\item the value $\omega_m(\Bh_m)$ of the universal quantum
	weight system is given by Eq.~(\ref{MainMain});
	
	\item the limit  $\lim_{q\to 1}\omega_m(\Bh_m)$ exists and is given by 
	Eq.~(\ref{AlmostMain});
	\item this limit coincides with $W_m(N)=w_\gl(\Bs_m)$.
\end{itemize}

It happens that computations with the universal quantum weight system
are quite complicated. To overcome this problem we introduce just another
function on Hecke algebras, called the characteristic mapping.
It looks more natural, can be easily computed and is
governed by much more simple recurrence relations. Unfortunately,
its values normally do not have a limit as $q\to1$,
which does not allow one to use it directly to deduce formulas for
the averaged values. Instead, we compute the average value of the
quantum characteristic mapping, rewrite the answer in terms of
universal quantum weight system, and then take the limit as $q\to1$,
which yields the desired formula.

\section{Preliminaries}\label{sglN}

\subsection{The universal $\gl$-weight system and statement of the main theorem}

For a definition of the universal $\gl$-weight system on permutations,
see~\cite{KLUI,KL,ZY}. Here we present an alternative but equivalent way 
to define it using matrix notation, see~\cite{GS}. Let $U$ be an $N$-dimensional vector space over the field $\mathbb{F}$, which can be either the field  $\mathbb{C}$
of complex numbers or the field $\mathbb{C}(q)$ of rational functions 
in one variable~$q$. 
Given an $N\times N$  matrix $X \in Mat_{N \times N}(\mathbb{F}) \otimes U$ we use notation $$
X_i = \BI\otimes \BI \otimes \ldots  \otimes\underset{i}{X}
\otimes \BI\otimes\ldots  \otimes \BI \in Mat_{N \times N}(\mathbb{F})^{\otimes k}\otimes U,
$$
where $\BI$ is the identity $N\times N$ matrix, for the matrix~$X$ acting on the 
tensor factor number~$i$.

Let $\BS_m$ be the group of permutations of $m$ elements. For a permutation 
$\alpha \in \BS_m$, we define the permutation matrix $P_{\alpha} \in End( \mathbb{C}^N)^{\otimes m}$ in the basis of decomposable tensors by the following rule $$P_{\alpha}:v_1 \otimes \dots \otimes v_m \mapsto v_{\alpha^{-1}(1)} \otimes \dots \otimes  v_{\alpha^{-1}(m)}.$$
Relations in the universal enveloping algebra $U({\gl(N)})$ can be rewritten in the matrix form as
$$E_1E_2 - E_2E_1 = P_{(12)}E_2 - E_2P_{(12)},$$
where $E= ||E_i^j||_{i,j=1}^N \in Mat_{N\times N}(\mathbb{C})\otimes U({\gl(N)})$ is the matrix of the algebra generators.

For $k\le n$, denote by $Tr_{k,n}$ the 
operator taking the trace in the tensor components $k,k+1,\dots,n$.
\begin{defn}
    For a permutation $\alpha \in \BS_m$, we set
    the value of the ${\gl(N)}$-\emph{weight system} equal to $$w_{{\gl(N)}}(\alpha) = Tr_{1,m}(E_1 \dots E_m P_{\alpha^{-1}}).$$
\end{defn}
The value  $w_{\gl(N)}(\alpha)$ of the ${\gl(N)}$-weight system 
belongs to the center $ZU(\gl(N))$ of $U({\gl(N)})$  and, therefore, can be expressed as a polynomial in the Casimir elements 
$$C_{m,N}=C_m = Tr_{1,m}(E_1 \dots E_m P_{(m-1 \, m)}P_{(m-2 \, m-1)}\dots P_{(1 \, 2)})=Tr(E^m) \in ZU({\gl(N)}).$$
\begin{st}\cite{KL,ZY}  \label{ClassUnivProp}
For any permutation $\alpha\in \BS_m$ there exists a universal polynomial $w_{\frak{gl}}(\alpha)$ in the variables $N,C_1,C_2,\dots$ such that the substitution of specific $N$ and $C_i = C_{i,N}$ transforms it into the polynomial $w_{{\gl(N)}}(\alpha)$. 
\end{st}
\begin{defn}\cite{KL,ZY} The value of the {\it universal ${\frak{gl}}$-weight system} on a permutation $\alpha$ is the polynomial $w_{\frak{gl}}(\alpha)$.
\end{defn}
The Casimir elements are expressed in terms of one-part Schur functions~$S_i$ 
by the Schur substitution formula
$$1-Nv-\sum_{m=1}^{\infty}C_mv^{m+1} =\left( \frac{1- \frac{N+1}{2}v}{1-\frac{N-1}{2}v}\right)^N \times \frac{\sum\limits_{i=0}^{\infty}v^i\left(1- \frac{N-1}{2}v\right)^{-i}S_i}{\sum\limits_{i=0}^{\infty}v^i\left(1- \frac{N+1}{2}v\right)^{-i}S_i},$$
where~$N$ is interpreted as a formal variable.
Following~\cite{KLUI}, we call the corresponding expression
the \emph{Schur substitution}. It is related to the Harish-Chandra isomorphism. {\it The Harish-Chandra isomorphism} $HC:ZU{\gl(N)}\rightarrow \mathbb{C}[x_1,\dots,x_N]^{\BS_N}$,
which identifies the center $ZU{\gl(N)}$ of the universal enveloping algebra $U\gl(N)$ of~$\gl(N)$
with the ring $\mathbb{C}[x_1,\dots,x_N]^{\BS_N}$ of symmetric polynomials in $x_1,x_2,\dots,x_N$,
 is defined on the Casimir elements by the following Perelomov--Popov formula
$$1-Nv-\sum_{m=1}^{\infty}C_mv^{m+1}=\prod_{i=1}^{N}\frac{1-\left(x_i+\frac{N+1}{2}\right)v}{1-\left(x_i+\frac{N-1}{2}\right)v}.$$
In variables $x_1,\dots,x_N$, the functions $S_i$ become one-part Schur polynomials. From now on, we will use the same notation for the central elements of the algebra $U{\gl(N)}$ and for their images under the Harish-Chandra isomorphism. 

Let us introduce a notation for the average value of the universal $\frak{gl}$-weight system on permutations of~$m$ elements:
$$W_m = \frac{1}{m!} \sum_{\alpha \in \BS_m}w_{{\frak{gl}}}(\alpha)=w_\gl(\Bs_m).$$
We are ready to formulate the main theorem, which was conjectured in \cite{KLUI}.
Recall our notation $\nu=N-1$.

\begin{thm}\label{MainThm} \it The generating function for the polynomials $W_m$ has the following form
     $$ \sum_{m=0}^{\infty}\frac{W_m}{(m+\nu)!}t^{m} = \left( \frac{e^{t/2}-e^{-t/2}}{t} \right)^{-\nu}\sum_{k=0}^{\infty}\frac{S_k}{(k+\nu)!}t^k.$$
In particular, the polynomials~$W_m$ are linear combinations of the
one-part Schur polynomials $S_1,S_2,\dots$.     
\end{thm}

Note that the coefficients in the expansion of $\left( \frac{e^{t/2}-e^{-t/2}}{t} \right)^{-\nu}$ in powers of $t$ are polynomials in $N=\nu+1$. Substituting a specific value of~$N$ and setting $S_i = S_i(x_1,\dots,x_N)$,
 we obtain a formula for the average value of the $\frak{gl}(N)$-weight system.

\begin{rem}
	The equation in Theorem~\ref{MainThm} is equivalent to Eq.~(\ref{AlmostMain}),
$$
W_{m}(N) = \sum_{l=0}^m\frac{(\nu+m)_l}{l!}B_l^{( \nu)}\left(\frac\nu2\right)S_{m-l}.
$$
This is the form of the main equation we are going to prove.	
\end{rem}

 The proof of the theorem is given in Section~\ref{sHgl}.

\subsection{Quantum universal $\gl$-weight system}

The goal of the present section is to define the quantum universal weight system
by introducing axioms for it. Its existence will be proven later. 

 In the quantum case, the role of the group algebra of the symmetric group is played by the Hecke algebra of type $A$ a detailed description
 of which can be found in the survey paper \cite{OP}. 
 
 In formulas below, all tensor products are considered over the field of rational functions $\mathbb{C}(q).$
 \begin{defn} The \textit{Hecke algebra}  $\mathbb{H}_m$ (of type~$A$) is a unital associative  algebra over $\mathbb{C}(q)$, with unit $1_{\mathbb{H}_m}$, generated by the Artin generators $g_i$, $i = 1,...,m-1$, that satisfy the following relations:
$$g_ig_j = g_jg_i, \quad  |i-j| > 1, $$ 
$$ g_ig_{i+1}g_i = g_{i+1}g_ig_{i+1}, \quad 1 \leq i \leq m-2,$$ $$g_i^2 = 1_{\mathbb{H}_m} + (q - q^{-1}) g_i, \quad 1 \leq i \leq m-2.$$ 
\end{defn}
Treating $q$ as a formal parameter, the Hecke algebra $\mathbb{H}_m$ is isomorphic to the group algebra $\mathbb{C}(q)[\BS_m]$ of the symmetric group $\BS_m$. 
 Set $$C_{ij} = g_{j-1}g_{j-2}\dots g_{i}, \quad j,i \in \mathbb{Z}_{> 0}, \,i<j, \quad C_{ii}=1.$$
 %, \quad C_{ij} = g_{i}g_{i+1}\dots g_{j-1}, \quad j<i.
 The element $C_{ij}$ can be considered as an element of each Hecke algebra
 $\BH_m$, for $m=j,j+1,\dots$.
 The monomials $C_{i_1 1}C_{i_2 2}\dots C_{i_k k}$, such that $1 \leq k \leq m$, $1 \leq i_l \leq l$, form a basis for the Hecke algebra $\mathbb{H}_m$, see~\cite{J}.

\begin{lemma} \label{UsefulLemma}
    Let $U$ be a $\mathbb{C}(q)$-vector space and let $f_m : \mathbb{H}_m \rightarrow U$, $m=1,2,\dots$, be a sequence of linear mappings. Suppose that for each $m \in \mathbb{Z}_{>0}$, any
     $x \in \mathbb{H}_m $, and any~$i$ such that $1 \leq i \leq m-1$, there is $y \in \mathbb{H}_{m-1}$ such that $$f_m(g_i x) - f_m(xg_i) = f_{m-1}(y).$$  Under these assumptions, the sequence of mappings 
    $f_m, m\geq 1$, is uniquely determined by the values
    of these functions on the elements 
    $C_{m_1k_1} \dots C_{m_lk_l}$, where $m_1 < k_1 < m_2 < k_2 < \dots <m_l < k_l.$
\end{lemma}
\begin{pf} First, we prove the following claim. Each basis element $u \in \mathbb{H}_n$ can be reduced by a sequence of conjugations and expansions to 
	a linear combination of the elements $ C_{m_lk_l} C_{m_{l-1}k_{l-1}}\dots C_{m_1k_1}$, $m_1 < k_1 < m_2 < k_2 < \dots <m_l < k_l$. We allow only conjugations by elements of the algebra  $\mathbb{H}_{m-1}$. After decomposing an element into a sum, we allow the terms to be conjugated by different elements. Elements $C_{m_ik_i}$ with non-intersecting segments $m_i,m_i+1,\dots, k_i$ commute, therefore the expression  $ C_{m_lk_l} C_{m_{l-1}k_{l-1}}\dots C_{m_1k_1}$ is equal to the desired $C_{m_1k_1} \dots C_{m_lk_l}$.

 We prove the claim by induction on $m$. The case $m=1$ is trivial.
 Suppose $m>1$, note that each element of the basis $u = C_{i_1 1}C_{i_2 2}\dots C_{i_m m}$ either lies in $\mathbb{H}_{m-1}$ or is equal to $xg_{m-1}y$, where $x,y \in \mathbb{H}_{m-1}$. If $u \in \mathbb{H}_{m-1}$, 
 then we can use the induction hypothesis. If $u = xg_{m-1}y$, the elements $x$ and $y$ are products of generators, therefore, we can invert $x$ and conjugate $u$ by $x^{-1}$. By expanding $yx$ as a linear combination of basis elements, $yx = \sum\limits_{i}a_iu_i$, $u_i\in \mathbb{H}_{n-1}$, $a_i \in \mathbb{C}(q)$,  we obtain the expression $x^{-1}ux=\sum\limits_{i} a_i g_mu_i$. By the induction hypothesis, any element $u_i$ can be reduced to a linear combination of $ C_{m_lk_l} C_{m_{l-1}k_{l-1}}\dots C_{m_1k_1}$, $m_1 < k_1 < m_2 < k_2 < \dots <m_l < k_l \leq m-1$. The element $g_m$ does not affect the calculation, 
 since we only used conjugation with those elements from the algebra $\mathbb{H}_{m-2}$, which commute with $g_{m-1}$. It remains to note that $$g_m C_{m_lk_l} C_{m_{l-1}k_{l-1}}\dots C_{m_1k_1} = \left\{ 
     \begin{array}{l}
          C_{m_lm} C_{m_{l-1}k_{l-1}}\dots C_{m_1k_1}, \quad k_l=m-1\\
          C_{ m-1 m}C_{m_l k_l} \dots  C_{m_1k_1},\quad k_l<m-1
    \end{array} 
 \right. .$$ 
 We used conjugations only by elements $x \in \mathbb{H}_{m-1}$ that are representable as products of the generators. For any such $x$ and any $u \in \mathbb{H}_{m}$, there is $z \in \mathbb{H}_{m-1}$ such that $f_m(x^{-1}ux) = f_m(u)+f_{m-1}(z)$. Thus, the assertion of the Lemma is reduced to the claim by a simple induction on $m$.\hfill\rule{6.5pt}{6.5pt}
 \begin{rem}
     Note that the claim gives an algorithm for representing an element of the algebra $\mathbb{H}_{m}$ as a linear combination of the basic elements 
     of the form $ C_{m_lk_l} C_{m_{l-1}k_{l-1}}\dots C_{m_1k_1}$. 

     Indeed, let $u \in  \mathbb{H}_{m} \setminus \mathbb{H}_{m-1}$ be a basis element. The element $u$ is uniquely represented as a product $xg_{m-1}y$, where $x,y\in\BH_{m-1}$, and after conjugating it with $x^{-1}$, we reduce it to the form $\sum a_i g_{m-1}u_i$, where $u_i \in \mathbb{H}_{m-1}$ are
      basis elements, $a_i \in \mathbb{C}(q)$. We then express each $u_i \notin \mathbb{H}_{m-2}$ as $\tilde{x_{i}}g_{m-2}\tilde{y_{i}}$. The elements $\tilde{y_{i}}$ have the form $C_{j\,m-2}$ for some $j \leq m-2$, and the elements $\tilde{x_{i}}$ have the form $C_{i_11}C_{i_22}\dots C_{i_{m-2} m-2 }$, thus the elements $\tilde{x_{i}}$ are invertible and the elements $\tilde{x_{i}}^{-1}$ can be explicitly expressed through Artin generators.
 
% \textcolor{red}{Что такое $\tilde x_i,\tilde y_i$? Что мы требуем от этих элементов?}     

      Conjugating $\tilde{x_{i}}g_{m-2}\tilde{y_{i}}$ with $\tilde{x}_{i}^{-1}$, which commutes with $g_{m-1}$, we get a sum of elements of the form $g_{m-1}g_{m-2}\tilde{u}_{ij}$ or $g_{m-1}\tilde{u}_{ij}$ where $\tilde{u}_{ij} \in \mathbb{H}_{m-2}$. Iterating this procedure, we obtain sums of expressions of the form
      $g_{i_1}g_{i_2} \dots g_{i_l}$, $i_1 > i_2 > \dots> i_{l}$. The normal form $g_{i_1}g_{i_2} \dots g_{i_l}$ is equal to $ C_{m_lk_l} C_{m_{l-1}k_{l-1}}\dots C_{m_1k_1}$ for some indices $m_i,k_i$.
 \end{rem}

% Note that any monomial $C_{k,n+1}$ for $k\leq n$ is equivalent to $g_k\dots g_{n-1} g_n = $\\
% $g_k\dots g_{n-2} g_n g_{n-1} = \dots= g_n g_{n-1} \dots g_{k} = C_{n+1,k}$ (all equalities are satisfied up to a cyclic shift of the generators).
% Therefore, it is sufficient to verify that any element can be reduced to a linear combination of monomials $C_{k_1n_1} \dots C_{k_ln_l}$, $n_1 < k_1 < n_2 < k_2 < \dots <n_l < k_l$ by cyclic shift of generators. Note that the indices are reversed. In particular, $f_{n}$ is completely determined by $f_{n-1}$ and the values on such elements.
    
%     We prove the claim by induction on $n$. The case $n=1$ is trivial.
% Suppose $n>1$, note that each element of the basis $C_{i_1 1}C_{i_2 2}\dots C_{i_n n}$ is equal to 
% $$
% \left\{ 
%     \begin{array}{l}
%           C_{i_n n-1}C_{i_1 1}\dots C_{i_{n-1} n} \, g_n \in \mathbb{H}_{n-1}(q) g_n,\quad i_n<n\\
%         C_{i_1 1}\dots C_{i_{n-1} n} \in \mathbb{H}_{n-1}(q),\quad i_n=n
%     \end{array} 
% \right. $$ 
% up to cyclic shift of the generators. Using the induction hypothesis, we obtain the element $C_{k_1n_1} \dots C_{k_ln_l}$ or $C_{k_1n_1} \dots C_{k_ln_l}g_n$ which is 
% $$\left\{ 
%     \begin{array}{l}
%           C_{k_1n_1} \dots C_{n n_l}, \quad k_l=n-1\\
%           C_{k_1n_1} \dots C_{k_ln_l} C_{n  n-1},\quad k_l<n-1
%     \end{array} 
% \right. $$ up to the cyclic shift of the generators.
%Note that on the n-th step, we did not use the cyclic property for the n-1-st generator. 

\end{pf}

Before introducing the quantum ${\frak{gl}}$-weight system, 
we define the universal characteristic map, which makes
computations more handleable.
For the next definition, we need the concatenation product on Hecke algebras.
 The \emph{concatenation product} is the algebra homomorphism $\mathbf{m}:\mathbb{H}_m \times \mathbb{H}_n \rightarrow \mathbb{H}_{m+n},\quad i=1,\dots,m,\quad j=1,\dots,n$, defined on the generators by the formula
$$\mathbf{m}: (g_i,g_j) \mapsto g_i g_{j+m}.$$

\begin{defn} \label{UnivCharMapDef}
    The { \it universal characteristic map} is the sequence of linear maps 
    $\chi_m : \mathbb{H}_m \rightarrow \mathbb{C}(q)[p_1, \dots, p_m]$,
    $m=1,2,\dots$, completely determined by the following properties:
        \begin{enumerate}
        \item~\label{UCMpm} $\chi_{m+n}(\mathbf{m}(x, y))=\chi_m(x)\chi_n(y)$, where $x \in \mathbb{H}_{m}$, $y \in \mathbb{H}_{n}$;
        \item\label{UCMpb} $\chi_m(C_{1m}) = p_m$, for each $m>0$;
        \item\label{UCMpe} $\chi_{m}(g_ix) = \chi_{m}(xg_i)$,
        for all $x \in \mathbb{H}_{m}$ and $1 \leq i \leq m-1$.
    \end{enumerate}
\end{defn}
Properties~\ref{UCMpm} and~\ref{UCMpb} imply the equality 
$$\chi_{m}(C_{m_1k_1} \dots C_{m_lk_l}) = p_{k_1-m_1+1}\dots p_{k_l-m_l+1} \times (p_1)^{m - \sum\limits_{i=1}^{l}(k_i-m_i+1)},$$
 where $m_1 < k_1 < m_2 < k_2 < \dots <m_l < k_l.$ Therefore, by Lemma \ref{UsefulLemma}, properties~\ref{UCMpm}--\ref{UCMpe} uniquely determine the set of mappings. Consistency of the definition will be proven in the next part of the section. 

\begin{defn}
    The {\it quantum trace} on $\mathbb{H}_m$, $m\geq 1$ is the sequence of linear mappings $$\BT r_{m}: \mathbb{H}_m \rightarrow \mathbb{H}_{m-1}$$ defined on the basic element $C_{i_1 1}C_{i_2 2}\dots C_{i_m m}$ by the formula
    $$\BT r_{m}(C_{i_1 1}\dots C_{i_m m})=
\left\{ 
    \begin{array}{l}
        (C_{i_1 1}\dots C_{i_{m-1} m-1}) C_{i_m m-1} ,\quad i_m<m\\
        q^{-1}[N]_q(C_{i_1 1}\dots C_{i_{m-1} m-1}),\quad i_m=m
    \end{array} 
\right. .$$
\end{defn}

\begin{defn} \label{UnivWeightSystemDef}
    The {\it universal quantum ${\frak{gl}}$-weight system} is the sequence of linear mappings $\omega_{m}:\mathbb{H}_m \rightarrow \mathbb{C}(q)[q^{-N};C_1,C_2,\dots],$ $m=1,2,\dots,$ completely determined by the following properties:
    \begin{enumerate}
        \item\label{UQWpm} $\omega_{m+n}(\mathbf{m}(x, y))= \omega_m(x)\omega_n(y)$, where $x \in \mathbb{H}_{m}$, $y \in \mathbb{H}_{n}$;
        \item\label{UQWpb} $\omega_m(C_{1m}) = C_m$, for any $m>0$;
        \item\label{UQWpe} $\omega_{m}(g_ix - xg_i)=$ $$= \omega_{m-1}(\BT r_{m}(g_{m-1}^{-1} \dots g_{i+1}^{-1} x g_{i+1} \dots g_{m-1}-g_{m-1}^{-1} \dots g_{i}^{-1} x g_{i} \dots g_{m-1})),$$
        for any $x \in \mathbb{H}_{m}$ and $1 \leq i \leq m-1$.
    \end{enumerate}
\end{defn}
The universal quantum ${\frak{gl}}$-weight system is uniquely determined by properties \ref{UQWpm}--\ref{UQWpe} due to the relation
$$\omega_{m}(C_{m_1k_1} \dots C_{m_lk_l}) = C_{k_1-m_1+1}\dots C_{k_l-m_l+1} \times (C_1)^{m - \sum\limits_{i=1}^{l}(k_i-m_i+1)},$$ where $m_1 < k_1 < m_2 < k_2 < \dots <m_l < k_l$ and Lemma \ref{UsefulLemma}. We have not yet establish that the maps $\chi_m$ and $\omega_m$ exist. Below in this section we will construct a characteristic mapping and a quantum weight system in the Reflection Equation algebra. This construction will allow us to introduce universal properties of the maps $\chi_m$ and $\omega_m$ (Proposition \ref{UnivPropCharMap2} and Corollary \ref{CorUnivProp}) and to prove their existence (Proposition \ref{PropUnivCharMap} and Proposition \ref{PropUnivWeightMap}).

\begin{ex}
Each basis element $C_{i_1 1}C_{i_2 2}\dots C_{i_m m}$ of the Hecke algebra can be associated with a permutation. This permutation is obtained by replacing each generator $g_i$ with the transposition $(i,i+1)$. This defines a bijective correspondence between permutations and the basis elements of the Hecke algebra.
However, the composition of the correspondence with the quantum weight system does not determine the weight system, see examples below.
We now compute the values of the quantum weight system for certain elements of the Hecke algebra corresponding to specific permutations.

    \begin{enumerate}
        \item The value of the quantum weight system on the standard cycle $(123)$ is $\omega_3(g_2g_1) = C_3$.
        \item  The value on the cycle $(132)$ is computed as follows
        \begin{align*}
           & \omega_3(g_1g_2) = \omega_3(g_2g_1) + \omega_2 (\BT r_{3} (g_2^{-1} g_2 g_2) - \omega_2(\BT r_{3}(g_2^{-1} g_1^{-1} g_2 g_1 g_2))= \\
           &= C_3 + \omega_2(1) - q^{-1}[N]_q \omega_2(g_1) = C_3 + C_1^{2} - q^{-1}[N]_q C_2.
        \end{align*} 
        Here we used the braid relation $g_2^{-1} g_1^{-1} g_2 g_1 g_2 = g_2^{-1}  g_2 g_1 g_2^{-1} g_2 = g_1$.

        \item The value on the transposition $(13)$ is \begin{align*}    
       & \omega_3(g_1g_2g_1) = \omega_3(g_2g_1^2)+\omega_2(\BT r_{3}(g_2^{-1}g_2g_1g_2)) - \omega_2(\BT r_{3}(g_2^{-1}g_1^{-1}g_2g_1g_1g_2)) =\\
       &= C_1C_2 + (q-q^{-1})C_3 +C_2 - \omega_2(\BT r_{3}(g_1g_2^{-1}g_1g_2)) = C_1C_2+(q-q^{-1})C_3\end{align*}
       The braid relation was applied twice $g_2^{-1}g_1^{-1}g_2g_1g_1g_2 = g_1g_2^{-1}g_1^{-1}g_1g_1g_2 = g_1^2g_2g_1^{-1}$.

        \item The permutation $(1243)$ corresponds to the element $g_1g_3g_2$ of the Hecke algebra. The value of the quantum weight system on this element is computed as:
        \begin{align*}
            &\omega_4(g_1g_3g_2) =  \omega_4(g_3g_2g_1)+ \omega_3(\BT r_{4} (g_3^{-1}g_2^{-1}g_3g_2^2g_3 - g_3^{-1}g_2^{-1}g_{1}^{-1}g_3g_2g_1g_2g_3)) = \\
            &= C_4 +C_1C_2 - q^{-1}[N]_qC_3 
        \end{align*}
        The calculation relies on: $g_3^{-1}g_2^{-1}g_{1}^{-1}g_3g_2g_1g_2g_3 = g_3^{-1}g_2^{-1}g_3g_2g_1g_3 = g_2g_3^{-1}g_2^{-1}g_2g_1g_3 = g_2g_1.$

        \item For the permutation $(13)(24)$, the corresponding element of the Hecke algebra looks like $g_2g_1g_3g_2$. 
        \begin{align*}
            &\omega_4(g_2g_1g_3g_2) = \\ 
            &=\omega_4(g_1g_3 + (q-q^{-1})g_1g_3g_2) + \omega_3(\BT r_{4}(g_3^{-1}g_1g_3g_2g_3 - g_3^{-1}g_2^{-1}g_1g_3g_2^2g_3)) = \\ 
            &= \omega_4(g_1g_3) + (q-q^{-1})\omega_4(g_1g_3g_2) + \omega_3(g_1g_2) - \omega_3(g_2+\frac{1}{q^{2N}}g_2g_1) = \\
            &= C_2^2 + (q-q^{-1})(C_4+C_2C_1) +C_1^2 - q^{-1}[N]_qC_2-C_1C_2
        \end{align*}
        To calculate the value $\omega_3(\BT r_{4}(g_3^{-1}g_2^{-1}g_1g_3g_2^2g_3))$, we used the following chain of equalities $g_3^{-1}g_2^{-1}g_1g_3g_2^2g_3 = g_2g_3^{-1}g_2^{-1}g_1(1+(q-q^{-1})g_2)g_3 = g_2g_3^{-1}g_2^{-1}g_3g_1 + (q-q^{-1})g_2g_3^{-1}g_2^{-1}g_1g_2g_3 = g_2^2g_3^{-1}g_2^{-1}g_1 + (q-q^{-1})g_2g_1g_2g_3g_2^{-1}g_1^{-1}$ and the trace property $\BT r_4(g_3^{-1}) = \BT r_4(g_3) +(q-q^{-1}) q^{-1}[N]_q= q^{-2N}$.

        \item For the permutation $(143)$, the corresponding element of the Hecke algebra looks like $g_1g_3g_2g_1$.
        \begin{align*}
           & \omega_4(g_1g_3g_2g_1)= \\
           &=C_3C_1 + (q-q^{-1})C_4 + \omega_3(\BT r_{4}(g_3^{-1}g_2^{-1}g_3g_2g_1g_2g_3-g_3^{-1}g_2^{-1}g_1^{-1}g_3g_2g_1^{2}g_2g_3))= \\ 
           &=C_3C_1 + (q-q^{-1})C_4 + \omega_3(\BT r_{4}(g_2g_3^{-1}g_1g_2g_3-g_2g_1g_3^{-1}g_2^{-1}g_1g_2g_3)) = \\
           &= C_3C_1+(q-q^{-1})C_4
        \end{align*}
        The last equality follows from: $g_3^{-1}g_2^{-1}g_1g_2g_3 = g_1g_2g_3g_2^{-1}g_1^{-1}$ and $g_3^{-1}g_2g_3 = g_2g_3g_2^{-1}$. The coefficient of $C_3$ vanishes: $-q^{-1}[N]_q +1 - q^{-2N} = 0$.

        \item For the permutation $(1423)$, the corresponding element of the Hecke algebra looks like $g_2g_1g_3g_2g_1$.
        \begin{align*}
            &\omega_4(g_2g_1g_3g_2g_1)= \\ 
            &= \omega_4(g_3g_2g_1) + (q-q^{-1})\omega_4(g_1g_3g_2g_1) + \omega_3(g_1g_2g_1 - g_2(q^{-1}[N]_qg_1 + (q-q^{-1})g_1) = \\
            &= C_4 + C_1C_2 - q^{-1}[N]_qC_3+(q-q^{-1})C_3C_1 + (q-q^{-1})^2C_4
        \end{align*}
        This combines results from items 3 and 6.

        \item For the permutation $(14)(23)$, the corresponding element of the Hecke algebra looks like $g_1g_2g_1g_3g_2g_1$.
        Using recurrence relations, we obtain 
        \begin{align*}
            &\omega_4(g_1g_2g_1g_3g_2g_1) = \\ 
            &= \omega_4(g_2g_1g_3g_2) + (q-q^{-1})\omega_4(g_2g_1g_3g_2g_1) + \omega_3(g_1^2g_2g_1) - \omega_3(g_1g_2g_1^2) = \\ 
            &=  \omega_4(g_2g_1g_3g_2) + (q-q^{-1})\omega_4(g_2g_1g_3g_2g_1) + \omega_3(g_2g_1) - \omega_3(g_1g_2)
        \end{align*}
        Substituting from items 5 and 7:
        \begin{align*}
            &\omega_4(g_1g_2g_1g_3g_2g_1) = \\ 
            &= C_2^{2} + (q-q^{-1})(2C_4 + C_1C_2-q^{-1}[N]_qC_3) +(q-q^{-1})^2C_3C_1 + (q-q^{-1})^3C_4
        \end{align*}
    \end{enumerate}
    In particular, we see that the quantum weight system 
    can take different values on cyclically equivalent permutations,
    as well as the four-term relations do not hold even for diagrams with two chords:
    $$\omega_4(g_1g_3) -\omega_4(g_2g_1g_3g_2) \neq \omega_4(g_1g_2g_1g_3g_2g_1) - \omega_4(g_2g_1g_3g_2). $$ Despite the quantum corrections, the zeroth order terms in $q$ of the obtained expressions coincide with the classical $\gl$ weight system, although at the moment we do not know how to prove this assertion
    in full generality. \hfill\rule{6.5pt}{6.5pt}
\end{ex}

From now on, let us identify the ranges of the invariants $\chi_m$
and $\omega_m$ by expressing the variables $C_1,C_2,\dots$
in the variables $p_1,p_2,\dots$ by means of the following
equality between the generating functions for them:
\begin{equation}\label{e-Cp}
e^{\frac{t}{q-q^{-1}}} \left(q^{-1}[N]_q + \sum\limits_{m=1}^{\infty}\frac{p_m}{m!}\frac{(-t)^m}{(q-q^{-1})^m}\right) = q^{-1}[N]_q + \sum\limits_{m=1}^{\infty}\frac{C_m}{m!}t^m;
\end{equation}
the coefficients are rational functions of~$q$.
This relation is explained using the quantum analogue of the Perelomov-Popov formula (see \cite{GPS}).
It allows one to express the values of the maps $\chi_m$ and $\omega_m$ 
through each other using a formula of the form
$$\omega_m(x) =\frac{1}{(q-q^{-1})^m} \sum\limits_{I \subset \{1,\dots,m\}}(-1)^{|I|}\chi_{|I|}(\tilde{x}_I), \quad  x \in \mathbb{H}_m,$$
where the sum is taken over all subsets $I$ of the set $\{1,\dots,m\}$, symbol $|I|$ denotes the cardinality of $I$ and the element $\tilde{x}_I \in \mathbb{H}_{|I|}$ is obtained from $x$ by conjugation with some invertible element of the algebra $\mathbb{H}_n$ and taking all traces in positions $|I|+1,\dots,m$. Details can be found in Proposition \ref{PropWeightCharConnection}.

\subsection{$R$-matrix representations of Hecke algebras}

In this Section, following~\cite{GPS,GS}, we introduce the Reflection
equation algebra~$\BM(N)$ and two sequences of linear mappings
$$
\chi_{m,N}:\BH_m\to Z\BM(N),\qquad \omega_{m,N}:\BH_m\to Z\BM(N),
$$ 
called the characteristic map and the quantum $\mathfrak{gl}(N)$-weight system,
respectively. Both mappings take value in the center $Z\BM(N)$
of $\BM(N)$, which can be identified with the ring of polynomials
in either of the two set of variables, $p_1,p_2,\dots,p_N$,
$p_i=\chi_{m,N}(C_{1m})$, or $C_1,\dots,C_N$, $C_i=\omega_{m,N}(C_{1m})$.
These two sets of variables are related by Eq.~(\ref{e-Cp}).
Similarly to the classical case, the invariants $\chi_{m,N}$ and
$\omega_{m,N}$ can be viewed as specializations of the universal
characteristic map~$\chi_m$ and the universal quantum $\mathfrak{gl}$-weight
system~$\omega_m$, respectively.
To construct the quantum ${\gl(N)}$-weight system and characteristic map, we need $R$-matrix representations of Hecke algebras (see \cite{OP}) and the associated Reflection Equation algebras (see \cite{GPS}).

We denote by $R=R(N)$, for $N=1,2,\dots$, the 
\emph{Drinfeld--Jimbo $R$-matrix}, which is the
$N^2\times N^2$-matrix acting
on the vector space~$\BC(q)^N\otimes\BC(q)^N$ according to the formula
\begin{equation}
	R(N) = q \sum\limits_{i=1}^N E_{i}^i \otimes E_{i}^i + \sum\limits_{i \neq j} E_{i}^j \otimes E_{j}^i + (q-q^{-1}) \sum\limits_{i<j}E_{i}^i \otimes E_{j}^j,
\end{equation}
where $E_i^j$, $i,j=1,\dots,N$ are standard matrix units. The $R$-matrix acting on the tensor
product of the $i$~th and the $(i+1)$~st factors in a tensor power $V^{\otimes m}$,
$V\equiv\BC(q)^N$, 
will be denoted $R_{i,i+1}$, $i=1,\dots,m-1$.

Below, we will make use of the following properties of the $R$-matrix:
\begin{itemize}
	\item $R^2=I+(q-q^{-1})R$ (the \emph{Hecke relation});
	\item $R_{i-1,i}R_{i,i+1}R_{i-1,i}
	=R_{i,i+1}R_{i-1,i}R_{i,i+1}$ (the \emph{braid relation}).
\end{itemize} 

For a vector space $V=\BC(q)^N$, the mapping 
\begin{equation}  
	\rho : \mathbb{H}_m \rightarrow End(V^{\otimes m}), \qquad 
	\rho:g_i \mapsto R_{i,i+1},
\end{equation} as can be easily seen, defines a representation 
of the Hecke algebra $\mathbb{H}_m$, which we call
the \textit{R-matrix representation}.

\begin{defn} \label{DefRTr}
 The \textit{R-trace} is the mapping $Tr_R : Mat_{N \times N}(\mathbb{C}(q))\otimes V \rightarrow V$ defined by the following formula:
    $$Tr_R(X) = Tr(DX), \quad  X \in Mat_{N \times N}(\mathbb{C}(q))\otimes V,$$ 
    where $D = diag(q^{1-2N},q^{3-2N}, \dots , q^{-1})$.
\end{defn}
 A direct generalization of Definition \ref{DefRTr} is 
that of the \emph{multiple $R$-trace} for a more general matrix $X \in Mat_{N \times N}(\mathbb{C}(q))^{\otimes m} \otimes U$, $1\le k\le n\le m$,
 $$\langle X \rangle_{k,n} = Tr_{R}(X) = Tr_{k,n} (D_k \dots D_n X),$$
 here we also introduce a convenient notation for the $R$-trace. 
The $R$-trace has the cyclic property 
$$\langle XR_{i,i+1} \rangle_{1,m} = \langle R_{i,i+1}X \rangle_{1,m}$$ for any matrix $X \in Mat_{N \times N}(\mathbb{C})^{\otimes m} \otimes U$ and $1 \leq i \leq m-1$, and it is normalized, so that 
$$\langle R_{1,2}\rangle_{2,2}=\BI, \qquad\langle \BI \rangle_{1,1} = q^{-1}[N]_q.$$

% The \emph{Harish-Chandra homomorphism} maps symmetric polynomials of $N$ variables to the center of the Reflection Equation algebra.

 The {\it Reflection Equation algebra} $\BM(N)$
  is generated by $N^2$ generators $m_i^j$, their permutation relations being written in the matrix form as
\begin{equation}\label{REfirst}
    R_{1,2} M_1 R_{1,2} M_1 = M_1 R_{1,2} M_1 R_{1,2},
\end{equation} where $M = ||m_i^j||_{1 \leq i,j, \leq N}$.

For a matrix  $ K \in Mat_{N \times N}(\mathbb{C}(q)) \otimes \BM(N)$,
 let us denote 
 $$K_{\overline{i}} = R_{i-1,i} \dots R_{1,2} K_1 R_{1,2}^{-1} \dots R_{i-1,i}^{-1},$$
 so that
 $$
 K_{\bar{i}}=R_{i-1,i}K_{\overline{i-1}}R_{i-1,i}^{-1},\qquad i=2,3,\dots,m.
 $$
  There is a sequence of mappings $\chi_{m,N} : \mathbb{H}_m \rightarrow \BM(N)$
   from the Hecke algebras $\mathbb{H}_m$, $m=1,2,\dots$, to the Reflection Equation algebra $\BM(N)$. They are defined by the formula
   $$\chi_{m,N}(x) = \langle M_{\overline{1}} \dots M_{\overline{m}}\rho (x) \rangle_{1,m}, \quad  x \in \mathbb{H}_m.$$
 \begin{defn}
     \rm The mappings $\chi_{m,N}, \,\, m \in \mathbb{Z}_{\geq 1}$ are called \textit{characteristic mappings}.
 \end{defn}
 The value $\chi_{m,N}(x)$ of the characteristic mapping can be expressed as a polynomial in
 the elements~$p_1,p_2,\dots$ defined by
 
 \begin{equation}\label{e-qps}
 	p_m=\chi_{m,N}(C_{1m}) = \langle M_{\overline{1}} \dots M_{\overline{m}}R_{m-1,m}R_{m-2,m-1}\dots R_{1,2}\rangle_{1,m}.
 \end{equation}
 The elements $p_i$ are called {\it quantum power sums}. They lie in the center of the Reflection Equation algebra and, therefore, the value $\chi_{m,N}(x)$ of the characteristic mapping is a central element of the algebra $\BM(N)$ for any $x \in \mathbb{H}_m$.

 \begin{st} \label{PropUnivCharMap}
     The universal characteristic mapping $\chi_m$ is well-defined;  the substitution
of Eq.~(\ref{e-qps})
transforms it into the characteristic mapping $\chi_{m,N}$.
 \end{st}
 \begin{pf} In the case $N > m$, the elements $p_1,\dots,p_m$ of the algebra $\BM(N)$ are algebraically independent. We deduce that the characteristic mapping $\chi_{m,N}$ defines a linear mapping $\mathbb{H}_m \rightarrow \mathbb{C}(q)[p_1,p_2,\dots, p_m] \hookrightarrow \mathbb{C}(q)[p_1,p_2,\dots]$. We check that this mapping satisfies all the properties of Definition \ref{UnivCharMapDef}, which implies the existence of a universal characteristic mapping.
 
     Let us use several standard facts about the Reflection Equation algebra (see \cite{GPS}). The cyclic property of the $R$-trace and the relations in the algebra $\BM(N)$ imply the equality   
     $$\langle  M_{\overline{1}} \dots M_{\overline{m}}\rho (x)R_{i,i+1}\rangle_{1,m} = \langle  M_{\overline{1}} \dots M_{\overline{m}}R_{i,i+1} \rho (x)\rangle_{1, m}, \quad  1 \leq i \leq m-1.$$ 
     This is equivalent to property~\ref{UCMpm} of the definition of the universal characteristic map. 
Property~\ref{UCMpb} of the definition coincides with the definition of the element $p_m$. Property~\ref{UCMpe} follows from the formula 
     \begin{eqnarray*}
         & & \langle M_{\overline{1}} \dots M_{\overline{m}}\rho(x) f(R_{k,k+1},R_{k+1,k+2},\dots, R_{m-1,m})\rangle_{1, m} = \\ &=&  \langle M_{\overline{1}} \dots M_{\overline{k-1}}\rho(x)\rangle_{1,k-1} \times\langle M_{\overline{1}} \dots M_{\overline{m-k+1}}f(R_{1,2},R_{2,3},\dots, R_{m-k+1,m-k+2})\rangle_{1,m-k+1},\end{eqnarray*} where $1 \leq k \leq m$, $x \in \mathbb{H}_{k-1}$,
         which is valid  for any noncommutative polynomial $ f(z_{1},z_{2},\dots, z_{m-k+1})$. \hfill\rule{6.5pt}{6.5pt}
 \end{pf}     

\begin{st}[Universal property of the characteristic mapping] \label{UnivPropCharMap2}
     The universal characteristic mapping is the unique sequence of mappings from Hecke algebras $\mathbb{H}_m$, $m=1,2,\dots$, 
     to the ring of polynomials in infinitely many variables $\mathbb{C}(q)[p_1,p_2,\dots]$ such that for any $x \in \mathbb{H}_m$
      and any specific value of $N$, the polynomial $\chi_m(x)$ coincides with the element $\chi_{m,N}(x)$ after substitution of Eq.~(\ref{e-qps}).
\end{st}
\begin{pf}
             Define the homomorphism $f: \mathbb{C}(q)[p_1,p_2,\dots] \rightarrow Z\BM(N)$ from the ring of polynomials to the center of the Reflection Equation algebra by Eq.~(\ref{e-qps}). The sequences of maps $f \circ \chi_{m}$ and $\chi_{m,N}$ satisfy the assumptions of Lemma \ref{UsefulLemma} and coincide on elements of the form $C_{m_1k_1} \dots C_{m_lk_l}$, where $m_1 < k_1 < m_2 < k_2 < \dots <m_l < k_l$, and therefore, by the lemma, they coincide everywhere.  Therefore, the universal characteristic mapping satisfies the universal property.
             
             Conversely, if some mapping satisfies the universal property, then it coincides with the mapping $\chi_{m,N}$ for $N>m$, and hence with the mapping $\chi_m$.\hfill\rule{6.5pt}{6.5pt}
\end{pf}

 A simple calculation shows that the matrix $L = \frac{1-M}{q-q^{-1}}$ satisfies the \emph{modified Reflection Equation} \cite{GS}:
\begin{equation} \label{MREsecond}
    R_1L_1R_1L_1 -L_1R_1L_1R_1 = R_1 L_1 - L_1 R_1 .
\end{equation}
\begin{defn}{\normalfont \cite{GS} }\label{d-qws}
    The \emph{ quantum ${\frak{gl}}(N)$-weight system} 
    $\omega_{m,N}:\BH_m\to Z\BM(N)$
    is defined by the formula $$\omega_{m,N}(x) = \langle  L_{\overline{1}} \dots L_{\overline{m}}\rho (x)\rangle_{1,m}, \quad  x \in \mathbb{H}_m.$$
\end{defn}
The value of the quantum ${\frak{gl}}_{N}$-weight system can be written as a polynomial in the quantum Casimir elements. The \emph{quantum Casimir elements}
 are a direct generalization of the classical ones, 
 and are defined by
 $$C_m = \omega_{m,N}(C_{1m}).$$
They lie in the center $Z\BM(N)$ of the Reflection Equation algebra $\BM(N)$ and, therefore, the value $\omega_{m,N}(x)$ of the quantum ${\gl(N)}$-weight system is a central element of the algebra $\BM(N)$ for any $x \in \mathbb{H}_m$.
\begin{st} \label{PropUnivWeightMap}
    The universal quantum ${\frak{gl}}$-weight system 
    $\omega_m:\BH_m\to\BC(q)[q^{-N};C_1,C_2,\dots]$ is well defined and coincides with the quantum ${\gl(N)}$-weight system under the substitution of specific $N$ and the corresponding Casimir elements.
\end{st}
\begin{pf}
Let us check the defining relations of the definition of the universal ${\frak{gl}}$-weight system for the quantum ${\gl(N)}$-weight system.

    The recurrence relation for the quantum ${\frak{gl}}(N)$-weight system (property~\ref{UQWpe}) follows from the cyclic property of the trace and the following relation in the Reflection Equation algebra
    $$L_{\overline{1}}\dots L_{\overline{m}}R_i - R_i L_{\overline{1}}\dots L_{\overline{m}} = L_{\overline{1}}\dots L_{\overline{i}}L_{\overline{i+2}}\dots L_{\overline{m}} - L_{\overline{1}}\dots L_{\overline{i-1}}L_{\overline{i+1}}\dots L_{\overline{m}}.$$
    Properties~ \ref{UQWpm} and~ \ref{UQWpb} of the definition follow from the equality
         \begin{eqnarray*}
         & & \langle L_{\overline{1}} \dots L_{\overline{m}}\rho(x) f(R_{k,k+1},R_{k+1,k+2},\dots, R_{m-1,m})\rangle_{1,m} = \\ &=&  \langle L_{\overline{1}} \dots L_{\overline{k-1}}\rho(x)\rangle_{1,k-1}
          \times\langle L_{\overline{1}} \dots L_{\overline{m-k+1}}f(R_{1,2},R_{2,3},\dots, R_{m-k+1,m-k+2})\rangle_{1,m-k+1},\end{eqnarray*}
          where $1 \leq k \leq m$, $x \in \mathbb{H}_{k-1}$,
          which is valid for any polynomial $ f(z_{1},z_{2},\dots, z_{m-k+1})$
          in  noncommuting variables. 

       For $N>m$, the quantum Casimir elements $C_1,\dots, C_m$ are algebraically independent. Consequently, the quantum $\frak{gl}(N)$-weight system defines a mapping $$\chi_{m,N}:\mathbb{H}_{m} \rightarrow \mathbb{C}(q)[C_1,C_2,\dots, C_m] .$$ 
        By Remark 1 after Lemma \ref{UsefulLemma}, we can express the value of the resulting mapping in terms of quantum Casimir elements using only conjugations and traces in the Hecke algebra. The resulting expressions depend polynomially on the parameter $q^{-N}$. 
        For each element $x \in \mathbb{H}_m$, we fix such an expression $\omega_m(x)$. Thus, we have obtained a mapping $\omega_{m}:\mathbb{H}_m \rightarrow \mathbb{C}(q)[q^{-N}][C_1,C_2,\dots, C_m]$, $ x\mapsto \omega_m(x)$. This mapping may not satisfy the definition of a quantum weight system and may even be non-linear. However, substitution of a given~$N$, $N>m$, transforms $\omega_{m}$ into the mapping $\omega_{m,N}:\mathbb{H}_{m}  \rightarrow \mathbb{C}(q)[C_1,C_2,\dots,C_m]$, which satisfies the definition of a quantum weight system.

        Let us check that the sequence of mappings $\omega_{m}$
        constructed above is a quantum weight system in the sense of Definition \ref{UnivWeightSystemDef}.
        For each element $x \in \mathbb{H}_m$, we choose a set of polynomials $P(x,f) \in \mathbb{C}(q)[q^{-N}]$. The polynomial $P(x,f)$ is the coefficient of the monomial $f = f(C_1,\dots,C_n)$ in the polynomial $\omega_m(x)$. 
        
        Each of the defining relations~ \ref{UQWpm}-- \ref{UQWpe} of Definition \ref{UnivWeightSystemDef} imposes some polynomial conditions on the polynomials $P(x,f)$:
        \begin{itemize}
           \item The linearity of the constructed mapping is equivalent to the property $P(\alpha x+\beta y,f) - \alpha P( x,f) -\beta P(y,f)=0$, for all $x,y \in \mathbb{H}_{n}(q)$, $\alpha,\beta \in \mathbb{C}(q)$.
            \item  Relation~\ref{UQWpm} is equivalent to the property $$\sum\limits_{gh=f}P(x,g)P(y,h) - P(\mathbf{m}(x,y),f)=0,$$
             for all 
             $$x \in \mathbb{H}_{k}, \quad y \in \mathbb{H}_{m-k},\quad 1 \leq k \leq m-1.$$ 
             Note that there are only finitely many pairs of monomials $g$ and $h$ such that $gh=f$. 
            \item Relation~\ref{UQWpb} is equivalent to the property $P(C_{1m},C_m) - 1=0$.
            \item  Relation~\ref{UQWpe} is equivalent to the property 
            \begin{align*}            	
            P(g_ix,f) & - P(xg_i,f)- P(\BT r_{m}(g_{m-1}^{-1} \dots g_{i+1}^{-1} x g_{i+1} \dots g_{m-1}),f)\\
            & + P(\BT r_{m}(g_{m-1}^{-1} \dots g_{i}^{-1} x g_{i} \dots g_{m-1}),f)=0,\end{align*}
        for any $x \in \mathbb{H}_{m}$ and $1 \leq i \leq m-1$.
        \end{itemize}

        For given elements of the Hecke algebra and a given monomial $f$, we denote by $P(x)$ the polynomial in one variable on the left-hand side of any of the identities.
        We know that $P(q^{-N})=0$ for all $N>m$, which means that the corresponding polynomial~$P$ is identically equal to zero. This means that the constructed mapping satisfies the definition of the universal quantum ${\frak{gl}}$-weight system.
        
         The proof of the fact that the universal quantum ${\frak{gl}}$-weight system after substitution of the Casimir elements coincides with the quantum ${\gl(N)}$-weight system completely repeats the proof of Proposition \ref{UnivPropCharMap2}.
 \hfill\rule{6.5pt}{6.5pt}

\end{pf}
\begin{cor}[Universal property of the quantum $\frak{gl}$-weight system]\label{CorUnivProp} 
 The universal quantum $\frak{gl}$-weight system is the unique sequence of mappings from Hecke algebras $\mathbb{H}_m$, $m=1,2,\dots$, 
 to the ring of polynomials in infinitely many variables $\mathbb{C}(q)[q^{-N}][C_1,C_2,\dots]$ such that for any $x \in \mathbb{H}_m$ and any particular $N$ the polynomial $\omega_{m}(x)$ coincides with the element $\omega_{m,N}(x)$ after substituting the corresponding value of $q^{-N}$ and the 
 corresponding quantum Casimir elements $C_i$.
\end{cor}
%\begin{rem}
%    The statement can be strengthened by requiring that the condition be satisfied for any Hecke $R$-matrix of $GL(N)$-type for which the notion of $R$-trace is defined.
%\end{rem}
\begin{pf}
%    Note that the calculations in the proof of Proposition \ref{PropUnivWeightMap} do not depend on the choice of a Hecke $R$-matrix. Therefore, conditions 1-4 of the universal $\frak{gl}$-weight system definition are satisfied in any Reflection Equation algebra of the $GL(N)$-type. This remark allows us to use the lemma \ref{UsefulLemma} as in the proof of Proposition \ref{UnivPropCharMap2}.
    
    It remains to verify the uniqueness of the mapping that satisfies the universal property. Let $\omega_m$, $m=1,2,\dots$, be any sequence
    of mappings satisfying assumptions of Definition~\ref{UnivWeightSystemDef}. As we have already noted, the quantum ${\gl(N)}$-weight system 
     determines a mapping $\mathbb{H}_{m} \rightarrow \mathbb{C}(q)[C_1,C_2,\dots, C_m] \hookrightarrow \mathbb{C}(q)[C_1,C_2,\dots]$ which coincides with the result of substituting a specific value of $N$ to the mapping $\omega_m$. Therefore, the polynomials $P(x,f)$ defined in the proof of Theorem \ref{PropUnivWeightMap} coincide with the corresponding coefficients of the polynomial $\omega_m(x)$ at infinitely many points. We conclude that the mappings $\omega_m$ coincide with the quantum weight system defined in the proof of Theorem \ref{PropUnivWeightMap}.\textbf{}\hfill\rule{6.5pt}{6.5pt}
\end{pf}

The corollary shows that all algebraic relations for the universal quantum $\frak{gl}$-weight system can be verified for the quantum $\frak{gl}(N)$-weight system in the Reflection Equation algebra. We see that 
the invariants $\chi_{m,N}$ and $\omega_{m,N}$ take values in the same ring, 
namely, the center $Z\BM(N)$ of the algebra $\BM(N)$. For generators of this ring we can use either elements $p_i$, $i=1,2,\dots$ or elements $C_i$, $i=1,2,\dots$. The connection between these generators is described by the following 
quantum analogue of the Harish-Chandra isomorphism and the Perelomov--Popov formula for this isomorphism.

The {\it Harish-Chandra isomorphism} for the Reflection Equation algebra is defined on the quantum power sums $p_n$  by the formula
\begin{equation} \label{HarishChandra}
    1+(q-q^{-1})\sum_{m=1}^{\infty}p_m v^m = \prod\limits_{i = 1}^{N} \frac{1 - q^{-2} \xi_i  v}{1 - \xi_i v}.
\end{equation}
In the variables $x_i=\frac{q^{1-N}-\xi_i}{q^2-1}$, the generating function of the Casimir elements $C_n$ has the form  
\begin{equation} \label{PerelomovPopov}1 - [N]_qv-(v-1+q^{-2})\sum_{m=1}^{\infty}q^{m+1}C_mv^{m} = \prod_{i=1}^N \frac{1-(x_i+[\frac{N+1}{2}]_q)v}{1-(q^2x_i+[\frac{N-1}{2}]_q)v}.\end{equation}
\begin{rem}
    Here we have identified the ranges of the classical and quantum weight systems. They take values in the ring of symmetric polynomials in variables $x_i$. Thus, the limit of the quantum weight system value as $q\rightarrow 1$ can be taken coefficient-wise.

    We consider all invariants over the ring $\mathbb{C}(q)[q^{-N},N]$, where $q$, $N$ and $q^{-N}$ are independent variables. If we want to pass to the limit, we should consider $q^{-N}$ as a power series in the neighborhood of $q= 1$ with coefficients that depend polynomially on $N$.
\end{rem}

% \begin{rem}
%     The invariant $h$ is the multiplicative analogue of the invariant $\omega$. The variables $\xi$ play the same role with respect to the variables $x$.
% \end{rem}

\begin{st} \label{PropWeightCharConnection}
    The value of the quantum weight system $\omega_{m,N}(x)$ on an element $x\in \mathbb{H}_m$ is expressed through the values of the characteristic mapping $\chi_{m,N}$ according to the formula
$$\omega_{m,N}(x) =\frac{1}{(q-q^{-1})^m} \sum\limits_{k=0}^{m} \left(\sum\limits_{1 \leq i_1 < i_2 < \dots <i_k \leq m} (-1)^{k}\chi_{k,N}(\tilde{x}_{i_1 i_2 \dots i_k})\right), \quad  x \in \mathbb{H}_m, $$
where the sum is taken over all directed subsets $1 \leq i_1 < i_2 < \dots <i_k \leq m$ of the set $\{1,\dots,m\}$. Denote 
the monomial $C_{1i_1}C_{2i_2}\dots C_{ki_k}$ by $r_{i_1 i_2 \dots i_k}$, for the set of indexes $ 1 \leq i_1 < i_2 < \dots <i_k \leq m$. The element $\tilde{x}_{i_1 i_2 \dots i_k} \in \mathbb{H}_{k}$ is obtained from $x$ using the formula $$\tilde{x}_{i_1 i_2 \dots i_k} = \langle r_{i_1 i_2 \dots i_k}^{-1}\,x \,r_{i_1 i_2 \dots i_k} \rangle_{i+1,m}.$$
\end{st}
\begin{pf}
We obviously have \begin{eqnarray*}
    \omega_{m,N}(x) = \langle L_{\overline{1}} \dots L_{\overline{m}} \rho(x)\rangle_{1,m} = \frac{1}{(q-q^{-1})^m}\langle (1-M_{\overline{1}}) \dots (1-M_{\overline{m}}) \rho(x)\rangle_{1,m} = \\ = \frac{1}{(q-q^{-1})^m}\sum\limits_{k=0}^{m} \left(\sum\limits_{1 \leq i_1 < i_2 < \dots <i_k \leq m} (-1)^{k}\langle M_{\overline{i_1}} \dots M_{\overline{i_k}} \rho(x)\rangle_{1,m}\right).
\end{eqnarray*}
    Note that $M_{\overline{k}} = R_{k-1,k}R_{k-2,k-1}\dots R_{l,l+1} M_{\overline{l}} (R_{k-1,k}R_{k-2,k-1}\dots R_{l,l+1})^{-1}$ for $l < k$.  Thus, we have $$M_{\overline{i_1}} \dots M_{\overline{i_k}} = \rho(C_{1i_1})M_{\overline{1}}\rho(C_{1i_1})^{-1}M_{\overline{i_2}} \dots M_{\overline{i_k}} = \rho(C_{1i_1})M_{\overline{1}}M_{\overline{i_2}} \dots M_{\overline{i_k}}\rho(C_{1i_1})^{-1},$$
    here we have used the mutual commutativity of $\rho(C_{ij})$ and $M_{\overline{k}}$ for $i,j<k$ or $i,j > k$. Repeating this procedure for all $M_{\overline{i_l}}$, we obtain
    $$M_{\overline{i_1}} \dots M_{\overline{i_k}} = \rho(r_{i_1 i_2  \dots i_k}) M_{\overline{1}} \dots M_{\overline{k}} \rho(r_{i_1 i_2  \dots i_k}^{-1}).$$ Application of the cyclic property of the $R$-trace to this formula proves the result.  \hfill\rule{6.5pt}{6.5pt}
\end{pf}
The formula is simplified for central elements of the algebra $\mathbb{H}_m$.
\begin{cor} \label{CorWeightCharConnection}
    For any central element $x \in \mathbb{H}_m$, we have $$\omega_{m,N}(x) =\frac{1}{(q-q^{-1})^m} \sum\limits_{k=0}^{m} (-1)^{k}\binom{m}{k}\chi_{k,N}\langle x\rangle_{k+1,m}.$$
\end{cor}

\section{Average value of the $\frak{gl}$-weight system}~\label{sqglN}
The average value of the ${\frak{gl}}(N)$-weight system over all permutations 
in the group $\BS_m$ can be expressed as the image of the group algebra $\mathbb{C}(q)[\BS_m]$ idempotent,
\begin{eqnarray*}
W_{m}(N) 
 &=& \frac{1}{m!} \sum_{\alpha \in \BS_m}w_{{\gl(N)}}(\alpha)\\
 &=&\frac{1}{m!}\langle E_1 \dots E_m \sum_{\alpha \in \BS_m}P_{\alpha^{-1}}\rangle_{1,m}\\
  &=& \frac{1}{m!} \langle E_1 \dots E_m \sum_{\alpha \in \BS_m}P_{\alpha}\rangle_{1,m}.
\end{eqnarray*}

\begin{defn} The {\it symmetrizer} in~$\BC[\BS_m]$ is the element $\Bs_m=\frac{1}{ m!}\sum_{\alpha \in \BS_m}\alpha$, so that $W_m(N)= w_{\gl(N)}(\Bs_m)$.
%An idempotent $\frac{1}{n!}\sum\limits_{\alpha \in S_n}\alpha \in \mathbb{C}[S_n]$ is called a symmetrizer. Its matrix image is denoted by $S^{(n)} = \frac{1}{n!}\sum_{\alpha \in S_n}P_{\alpha}$.
\end{defn}
The symmetrizer is a central idempotent of the algebra $\mathbb{C}[\BS_m]$.
It can be defined recursively by
$$
\Bs_m=\frac1m\Bs_{m-1}(id_m+(m-1)\cdot(m-1,m))\Bs_{m-1},
$$
where $id_m\in\BS_m$ is the identity permutation and $(m-1,m)\in\BS_m$
is a transposition.
Hence, its matrix image $S^{(m)} = \frac{1}{ m!}\sum_{\alpha \in \BS_m}P_{\alpha}
\in End(\BC^N)^{\otimes m}$ also 
can be calculated recursively by the following rule
$$S^{(1)} =1, \quad S^{(m)} = \frac{1}{m} S^{(m-1)}(1+(m-1)P_{(m-1,m)})S^{(m-1)}, \quad m \geq 2.$$
The corresponding central idempotent $\Bh_m$ in the Hecke algebra $\BH_m$ is called the $q$-\emph{symmetrizer} and is constructed similarly: 
$$\Bh_1 =1, \quad \Bh_m = \frac{1}{[m]_q} \Bh_{m-1}(q^{2-2m}+q^{-1}[m-1]_q g_{m-1})\Bh_{m-1}, \quad m \geq 2.$$
Let $H^{(m)} = \rho(\Bh_m)$ denote the $R$-matrix image of the idempotent $\Bh_m$.
\begin{defn} {\it The average value of the quantum ${{\gl(N)}}$-weight system} is the element 
	$$\Omega_m(N) = \omega_{m,N}(\Bh_m) = \langle L_{\overline{1}} \dots L_{\overline{m}} H^{(m)}\rangle_{1, m}.$$
\end{defn}

\begin{thm}\label{LimitThm}
    The image of $W_{m}(N)$ under the Harish-Chandra homomorphism coincides with the image of $\Omega_{m}(N)$ in the limit $q \rightarrow 1$.
\end{thm}

Along with the average value $W_m(N)$ of the ${\mathfrak gl}(N)$-weight system
and its quantum analogue $\Omega_m(N)$, we introduce one more pair of elements
$\widehat W_m(N)\in ZU\gl(N)$ and $\widehat \Omega_m(N)\in Z\BM(N)$. %(see \cite{OK}, \cite{JLM}).
The elements $W_m(N),\widehat W_m(N),\Omega_m(N),\widehat \Omega_m(N)$ 
can be considered as symmetric polynomials
in the set of independent variables $x_1,\dots,x_N$. Now, we have
\begin{itemize}
	\item there are known explicit formulas for $\widehat W_m(N), \widehat \Omega_m(N)$, which
	imply, in particular, that 
	$$\widehat W_m(N)=\lim_{q\to1}\widehat \Omega_m(N);$$
	\item the two sequences of polynomials $\Omega_m(N)$, $\widehat \Omega_m(N)$, $m=1,2,\dots$,
	are related by a triangular invertible linear transformation.
	Each coefficient of the transformation has a limit as $q\to1$
	and the limit transformation (which is also linear and triangular) 
	relates the corresponding sequences of polynomials $W_m(N)$ and $\widehat W_m(N)$.
	\end{itemize} 
	
	By combining the two assertions above, we establish the assertion of 
	Theorem~\ref{MainThm}.

In more detail, we set %(Proposition 3.1 \cite{JLM})
$$\widehat \Omega_{m}(N)= \langle L_{\overline{1}} (L_{\overline{2}}-q^{-1}[1]_q) \dots (L_{\overline{n}}-q^{-1}[m-1]_q)H^{(m)}\rangle_{1, m}\in Z\BM(N).$$ 
 The image of this element under the Harish-Chandra isomorphism is equal to %(Theorem 3.2 \cite{JLM})
\begin{eqnarray}\nonumber
    \frac{1}{(q^2-1)^m}\sum\limits_{1 \leq i_1 \leq \dots \leq i_m \leq N} \prod\limits_{k=1}^{m}(q^{-2(\nu-i_k+k)}-\xi_{i_k}) = \\ =\sum\limits_{1 \leq i_1 \leq \dots \leq i_m \leq N} \prod\limits_{k=1}^{m}  \left(x_{i_k}-q^{-N-1}\left[\frac{\nu}{2}-i_k+k\right]_q\right). \label{QQImm}
\end{eqnarray}
The elements $\widehat \Omega_{m}(N)$ were introduced in \cite{JLM}. Here, we adopt the definition given in Proposition 3.1 of \cite{JLM}, and Eq.~(\ref{QQImm}) restates Theorem 3.2 therein.

Similarly, set %(formula (2.1) \cite{OK})
 $$\widehat W_m(N)= \langle E_1 (E_2-1)\dots (E_m-m+1)S^{(m)}\rangle_{1,m} \in U\gl(N).$$ 
 The image of this element under the Harish-Chandra isomorphism is equal to %(formula (2.2) \cite{OK})
\begin{equation}\label{ClassImm}
    \sum\limits_{1 \leq i_1 \leq \dots \leq i_n \leq N} \prod\limits_{k=1}^{m}\left(x_{i_k}-\frac{N+1}{2}+i_k-k+1\right).
\end{equation}
The elements $\widehat W_m(N)$ are described in \cite{OK} (2.1)-(2.2).

    Passing to the limit $q\rightarrow1$ in Eq.~(\ref{QQImm}), we obtain Eq.~(\ref{ClassImm}). To prove Theorem \ref{LimitThm}, we need a couple of lemmas.
\begin{lemma}
    The following formulas hold
        \begin{equation}\label{TrH}       
    \langle H^{(m)}\rangle_{m,m} = q^{-1} \frac{[\nu+m]_q}{[m]_q} H^{(m-1)}, \quad 
    Tr_{m,m} (S^{(m)} )= \frac{\nu+m}{m} S^{(m-1)}.\end{equation}
\end{lemma}
\begin{pf}
    Using the definition of the $q$-symmetrizer, we obtain $$\langle H^{(m)} \rangle_{m,m} = \frac{q^{1-m}}{[m]_q}H^{(m-1)} (q^{1-m}\langle \BI^{\otimes m} \rangle_{m,m} + q^{m-2}[m-1]_q \langle R_{m-1,m} \rangle_{m,m})H^{(m-1)}. $$
    Since $\langle \BI^{\otimes m} \rangle_{m,m} = q^{-1}[N]_q\BI^{\otimes m-1}$ and $\langle R_{m-1,m} \rangle_{m,m}= \BI^{\otimes m-1}$, it follows that $$\langle H^{(m)} \rangle_{m,m} = (H^{(m-1)})^2 \frac{q^{1-m}(1-q^{-2N})+q^{m-1}-q^{1-m}}{q^m - q^{-m}} = q^{-1}\frac{1-q^{-2(\nu+m)}}{1-q^{-2m}}H^{(m-1)}.$$
    In the last implication we used the fact that $\Bh_{m-1}$ is a Hecke algebra idempotent.
    Passing to the limit $q\rightarrow 1$, we obtain the relation in the classical case. \hfill\rule{6.5pt}{6.5pt}
\end{pf}   

\begin{st} \label{PropImmanantWConnection}
    The elements $\widehat \Omega_m(N)$ admit the following expression in terms of $\Omega_m(N)$:
    \begin{equation} \label{qImm->WS}
        \widehat \Omega_m(N) = \sum\limits_{k=0}^m(-1)^k\frac{([\nu+m]_q)_k}{([m]_q)_k} e_k \Omega_{m-k}(N),
    \end{equation}
    where $e_k = q^{-k}\sum\limits_{1 \leq i_1<\dots < i_k \leq m}\prod\limits_{l=1}^{k}[i_l]_q$.
    A similar formula holds for the elements $\widehat W_m(N)$ and $W_m(N)$ in the algebra $U({\gl(N)})$:
     \begin{equation} \label{Imm->WS}
        \widehat W_m(N) = \sum\limits_{k=0}^m(-1)^k e_k \frac{(\nu+m)_k}{(m)_k} W_{m-k}(N),
    \end{equation}
    where $e_k = \sum\limits_{1 \leq i_1<\dots < i_k \leq m}\prod\limits_{l=1}^{k}i_l$.
\end{st}
\begin{pf}
   As in Proposition~\ref{PropWeightCharConnection}, any set of matrices $L_{\overline{i_1}} \dots L_{\overline{i_k}}$ is conjugate to a set of matrices $L_{\overline{1}} \dots L_{\overline{k}}$, and, therefore, the cyclic property of the $R$-trace and centrality of $H^{(m)}$ imply the equality $$\langle L_{\overline{1}} (L_{\overline{2}}-q^{-1}[1]_q) \dots (L_{\overline{m}}-q^{-1}[m-1]_q)H^{(m)}\rangle_{1, m}=\sum\limits_{k=0}^{m} (-1)^k e_k\langle L_{\overline{1}} \dots L_{\overline{m-k}}H^{(m)}\rangle_{1, m}.$$
    Using Eq.~(\ref{TrH}), we obtain the required result.
    Calculation in the case of the algebra $U({\gl(N)})$ is completely analogous.\hfill\rule{6.5pt}{6.5pt}
\end{pf}
The coefficients in formula (\ref{Imm->WS}) are the limits of the corresponding coefficients in formula (\ref{qImm->WS}). This remark is the key in the proof of Theorem \ref{LimitThm}.

\noindent
\textbf{Proof of Theorem \ref{LimitThm}} The transition matrix from
the elements $\widehat \Omega_m(N)$ to $\Omega_m(N)$ is unitriangular, which means that the elements of the inverse matrix are polynomials with integer coefficients in the elements of the original matrix. Moreover, the limits of the elements of the transition matrix in the quantum case are the corresponding elements in the classical case. Rewriting the elements $\Omega_m(N)$ as linear combinations of $\widehat \Omega_m$'s, we deduce that the limits of $\Omega_m(N)$ are equal to $W_m(N)$.\hfill\rule{6.5pt}{6.5pt}

%\begin{rem}
%    The element $W^{n,N}_q$ is a specific value of the universal polynomial $W^n_q = \omega_{n}(S_n)$ in variables $C_n$ and $q^{-N}$. Changing the variables $C_n$ with $C_n$ and passing to the limit $q \rightarrow 1$, we obtain the polynomial $W^n$. Expression $q^{N}$ should be treated as series arond $q=1$ with polynomial coefficients in $N$.
%\end{rem}

\section{Calculation of the quantum and classical ${\frak{gl}}$-weight systems average values}~\label{sHgl}
In this section we will calculate the average value of the quantum ${\frak{gl}}$-weight system and its limit. For any fixed $N\geq 1$, the value of the averaged quantum ${\gl(N)}$-weight system is calculated by the formula 
 \begin{eqnarray}\nonumber
        \Omega_m(N) &=& \omega_{m,N}(\Bh_m) =\frac{m!}{ (\nu+m)!}\frac{[\nu+m]_q!}{[m]_q!} \sum_{l=0}^{n} (\nu+m)_l \times \\ &&\times\frac{1}{q^{2l} l!}\beta_{l}^{(\nu+m-l,\nu)}\left(\frac{\nu}{2}\right) S_{m-l}(x_1,\dots,x_N) .\label{AvgVal}
    \end{eqnarray}
    Recall that $\nu=N-1$ and $\beta_{m}^{(h,k)}(\xi)$ is the q-Bernoulli polynomial of order $k$,
       $$\beta_{m}^{(h,k)}(\xi) = \frac{1}{(1-q^{-2})^m} \sum\limits_{i=0}^{m}\binom{m}{i}(-q^{-2\xi})^i \frac{(i+h)_k}{( [i+h]_q)_k }.$$ 
    Here $\xi$ is a formal parameter which can be specialized at integer and half-integer points, $m \in \mathbb{Z}_{\geq 0}$, $h,k \in \mathbb{Z}_{\geq 1}$ and $h \geq k$.

The limit of the $q$-Bernoulli numbers $\beta_{m}^{(1,1)}= \beta_{m}^{(1,1)}(0)$ of order $1$ is well known and coincides with the classical Bernoulli numbers. The $q$-Bernoulli polynomial of order $k$ can be expressed in terms of $q$-Bernoulli numbers of order 1, which allows us to calculate the limit $$\lim\limits_{q\rightarrow 1}\beta_m^{(h,k)}(\xi)=B_m^{(k)}(\xi)$$
and prove that it does not depend on the value of $h$ and is given by the generating series
$$\sum_{m=0}^\infty B_m^{(k)}(\xi)\frac{t^m}{m!}=e^{\xi t}\left(\frac{e^t-1}{t}\right)^{-k},$$
see Proposition~\ref{BkmLimit} below.

We conclude that all the summands in the right hand side of  Eq.~(\ref{AvgVal}) have limits as $q\to 1$ and the limit values satisfy the relation
$$W_m(N)=w_{\frak{gl}(N)}(\Bs_m)=\sum_{l=0}^m \frac{(\nu)_l}{l!}
B_l^{(\nu)}\left(\frac{\nu}{2}\right)S_{m-l}(x_1,\dots,x_N).
$$
Collecting these equalities for all $n$ in a generating series, we obtain the required relation of Theorem \ref{MainThm} for the ${\gl(N)}$-weight system. 

Proposition \ref{ClassUnivProp} and Corollary \ref{CorUnivProp} show that if a polynomial relation holds for all ${\gl(N)}$-weight systems, then it also holds for the universal ${\frak{gl}}$-weight system, so that we obtain expressions (\ref{AlmostMain}) and (\ref{MainMain}) for the classical and quantum universal ${\frak{gl}}$-weight systems, respectively. 

\subsection{Computation of the average value of the quantum $\gl(N)$-weight system}

In this section, we will prove that the average value of the quantum ${\gl(N)}$-weight system is given by Eq.~(\ref{AvgVal}).
We start with a formula from~\cite{GPS} for the mean value of the quantum 
character.

\begin{st}[\cite{GPS}] \label{REAcenter} In the Reflection Equation algebra, the following relation holds
    $$\chi_m(\Bh_m)=\langle M_{\overline{1}} \dots M_{\overline{m}}H^{(m)}\rangle_{1, m}  = q^{-m}S_{m}(\xi_1,\dots,\xi_N) = q^{-m} \sum\limits_{1 \leq i_1 \leq\dots \leq i_m \leq N}\xi_{i_1}\dots\xi_{i_m}.$$

\end{st}
%In this formula, we have identified a central element 
%of the Reflection Equation algebra
%with its image under the Harish-Chandra isomorphism. 

\begin{st}
    In terms of the $\xi$  variables, the value of the quantum 
    ${\frak{gl}}$-weight system on the $q$-symmetrizer is as follows
        \begin{equation}\label{Quant}
        \omega_{m,N}(\Bh_m) = \Omega_m(N) = \frac{1}{(q^2-1)^m}\sum\limits_{k=0}^{m}(-1)^k  \binom{m}{k}\frac{([\nu+m+l]_q)_{m-k}}{([m+l]_q)_{m-k}}S_k(\xi_1,\dots,\xi_N).
    \end{equation}
\end{st}
\begin{pf} Combining Corollary \ref{CorWeightCharConnection} and formula (\ref{TrH}) we obtain
$$ \Omega_m(N) = \frac{1}{(q^2-1)^m}\sum\limits_{k=0}^{m}(-1)^k  \binom{m}{k}\frac{([\nu+m+l]_q)_{m-k}}{([m+l]_q)_{m-k}}\chi_{k,N}(\Bh_k).$$
It remains to apply Proposition \ref{REAcenter} to the obtained expression.
\hfill\rule{6.5pt}{6.5pt}
    %We used the same calculation as in the proof of Proposition \ref{PropImmanantWConnection}. Here, $L_{\overline{k}} = \frac{1-M_{\overline{k}}}{q-q^{-1}}$ and $e_k$ are replaced by $\binom{n}{k}$.\hfill\rule{6.5pt}{6.5pt}
\end{pf}

\begin{st}
    The one-part Schur polynomials have the property 
    \begin{equation}\label{Class}
        S_m(\xi_1,\dots,\xi_N)= \sum\limits_{k=0}^{m}(-1)^k \binom{\nu+m}{m-k} S_k(u - \xi_1,\dots,u - \xi_N)u^{m-k}.
    \end{equation}
\end{st}
% \begin{pf}
%     Let $X = diag(u-\xi_1, \dots,u-\xi_N)$ be a matrix of commutative variables. Note that $$s_n(u-\xi_1,\dots,u-\xi_N)= Tr_{(1\dots n)}(X_1 \dots X_n S^{(n)}).$$ A simple calculation with the expression $$s_n(\xi_1,\dots,\xi_N) = Tr_{(1\dots n)}((u-X_1) \dots (u-X_n ) S^{(n)})$$
%  shows that $$s_n(\xi_1,\dots,\xi_N) =\sum\limits_{k=0}^{n}u^{n-k}(-1)^{k}\binom{n}{k}Tr_{(1\dots n)}(X_1 \dots X_k S^{(n)}).$$
%  If we combine this with with formula (\ref{TrH}), we get formula (\ref{Class}).\hfill\rule{6.5pt}{6.5pt}
% \end{pf}
\begin{pf} This results from the following simple calculation with generating series:
    \begin{eqnarray*}
\sum_{m=0}^\infty S_m(\xi_1,\dots,\xi_N)t^m&=&\prod_{i=1}^N\frac{1}{1-\xi_i t}=\prod_{i=1}^N\frac{1}{1-ut+(u-\xi_i) t}
\\
&=&(1-u t)^{-N}\prod_{i=1}^N\frac1{1+(u-\xi_i) \frac{t}{1-u t}}\\
&=&
 (1-u t)^{-N}\sum_{m=0}^\infty S_m(u-\xi_1,\dots,u-\xi_N)\left(\frac{-t}{1-u t}\right)^m.
     \end{eqnarray*}
     Expanding $ (1-u t)^{-N}\left(\frac{-t}{1-u t}\right)^m$ into a power series in~$t$, we obtain the desired formula.\hfill\rule{6.5pt}{6.5pt}
\end{pf}
\begin{cor} \label{CorWFormula}
    Equation~(\ref{AvgVal}) is valid.
       
\end{cor}
\begin{pf} The coefficient of $S_{l}$ in the expansion of $\Omega_{n}(N)$ is the sum of the expressions
    %The coefficient of $s_{l}$ is the sum of the products of the coefficients from the expansions of $s_{k}$, $k\geq l$ by the formula (\ref{Class}) and the coefficients of $s_{k}$ in the formula (\ref{Quant}), that is, the sum of the following expressions
    $$(-u)^{k-l}\frac{1}{(q^2-1)^{m-l}}\frac{m!}{l!(k-l)!(m-k)!} \times \frac{([\nu+m]_q)_k (\nu+k)_{k-l}}{([m]_q)_{m-k} (k)_{k-l}}.$$
    This expression is obtained by substituting Eq.~(\ref{Class}) into 
    Eq.~(\ref{Quant}).
    After dividing the coefficient by $\frac{m!}{[m]_q!}\frac{[\nu+m]_q!}{(\nu+m)!} \frac{(\nu+m)_l}{(m-l)!}$, we get
        $$(-u)^k\frac{1}{(q^2-1)^{m-l}}\frac{[k]_q!}{[\nu+k]_q!} \frac{(\nu+k)!}{k!} \frac{(m-l)!}{(m-k)!(k-l)!}.$$
        The sum of these expressions for $k = l,l+1,l+2,\dots, m$ coincides with the $q$-Bernoulli polynomial. Note that by the definition of the variables $x_i$ we need to choose $u=q^{-\nu} = q^{-2\cdot\frac{\nu}{2}}$. \hfill\rule{6.5pt}{6.5pt}
\end{pf} 

\subsection{The limit of the $q$-Bernoulli polynomials}

%In this part of the section we will focus on calculating the limit of $q$-Bernoulli polynomials as $q \rightarrow 1$. Denote the $m$~th $q$-Bernoulli number of order $k$ by $\beta_{m}^{(h,k)} = \beta_{m}^{(h,k)}(0)$.

This section is devoted to a proof of the following result.

\begin{st} \label{BkmLimit}
	The $q$-Bernoulli polynomial $\beta_m^{(h,k)}(\xi)$
	has a limit as $q\to1$. This limit
	$$
	B_m^{(k)}(\xi)=\lim_{q\to1}\beta_m^{(h,k)}(\xi)
	$$
is independent of~$h$ and coincides with the classical Bernoulli polynomial,
whose sequence is defined by the generating series
$$
\sum_{m=0}^\infty B_m^{(k)}(\xi)\frac{t^m}{m!}=e^{\xi t}\left(\frac{e^t-1}{t}\right)^{-k}.
$$ 
\end{st}	

The proof is based on the following recursions for the 
$q$-Bernoulli polynomials. Denote by $\beta_m^{(h,k)}=\beta_m^{(h,k)}(0)$,
the corresponding $m$~th $q$-Bernoulli number.

\begin{lemma} \label{recurrence}
	We have 
\begin{eqnarray}\label{Polynomial->Number}
\beta_m^{(h,k)}(\xi)&=&\sum_{i=0}^m {m\choose i} \beta_{m-i}^{(h+i,k)}
([\xi]_q)^i;\\ \label{k->1}
\beta_m^{(h,k)}&=&\sum_{i_1+\dots+i_k=m} {m\choose i_1,\dots,i_k}
\prod_{j=1}^k \beta_{i_j}^{(h-j+1+i_1+\dots+i_{j-1},1)};\\ \label{h->1}
\beta_m^{(h+1,1)}&=& \beta_{m}^{(h,1)}-(1-q^{-2})\beta_{m+1}^{(h,1)}.
\end{eqnarray}
	
\end{lemma}	

The first relation expresses the $q$-Bernoulli polynomials in terms of
$q$-Bernoulli numbers. The second one reduces computation of
$q$-Bernoulli numbers to the case $k=1$, and reduces the latter to the
case $h=1$. Thus, all these relations together express $q$-Bernoulli 
polynomials in terms of the numbers $\beta_{m}^{(1,1)}$.

The relations of the last Lemma have classical analogues.
Namely, denoting $B_m^{(k)}=B_m^{(k)}(0)$, $k=0,1,2,\dots$, we have
\begin{eqnarray} \label{ClassPol->Number}
B_m^{(k)}(\xi)=\sum_{i=0}^m{m\choose i} B_{m-i}^{(k)}\xi^i; \\ \label{Clask->1}
B_m^{(k)}=\sum_{i_1+\dots+i_k=m}{m\choose i_1,i_2,\dots,i_k}
\prod_{j=1}^k B_{i_j}^{(1)}.
\end{eqnarray}
The coefficients in these relations are obtained by taking the limit $q\to1$.
Finally, we use the following result:

\begin{st}[\cite{C}] \label{b11Limit}
   For any $m>0$ the limit as $q \rightarrow 1$ of the $q$-Bernoulli number $\beta_{m}^{(1,1)}$ of order $1$  exists and is equal to the classical Bernoulli number $B_m^{(1)}(0)$.
\end{st}

Applying this Proposition and induction we obtain, step by step, 
the following equations:
$$
\lim_{q\to1} \beta_m^{(h,1)}=B_m^{(1)};
$$
$$
\lim_{q\to1} \beta_m^{(h,k)}=B_m^{(k)};
$$
$$
\lim_{q\to1} \beta_m^{(h,k)}(\xi)=B_m^{(k)}(\xi).
$$

Thus, it remains to prove Lemma \ref{recurrence} to complete the proof of Proposition \ref{BkmLimit}.

Let $q^{-N}$ be a formal variable. Define the linear operator $\Phi_y:  \mathbb{C}(q)[q^{-2y}] \rightarrow \mathbb{C}(q)$ by the formula
$$\Phi_y(q^{-2ay}) = \frac{a+1}{[a+1]_q}, \quad a \in \mathbb{Z}_{\geq 0}.$$
% , \quad \Phi_y(q^{2y})=\frac{0}{[0]_q}=0
This operator is a model of the {\it $p$-adic $q$-Volkenborn integral} (see \cite{KCR} and references therein). 
\begin{st} \label{PropIntegralFormula}
    The $q$-Bernoulli polynomial admits the following description in
    terms of the operator $\Phi_y$:
         \begin{equation} \label{Bernoulli}
             \beta_{m}^{(h,k)}(\xi) = \Phi_{y_1}\dots\Phi_{y_k} \left( q^{-2\sum\limits_{i=1}^k(h-i)y_i}([\xi+y_1+\dots+y_k]_q)^m\right).
        \end{equation}
\end{st}
\begin{pf} A simple calculation shows that
     \begin{eqnarray*}
        q^{-2\sum\limits_{j=1}^k(h-j)y_j}([\xi+y_1+\dots+y_k]_q)^m= \frac{q^{-2\sum\limits_{j=1}^k(h-j)y_j}}{(1-q^{-2})^m}(1 - q^{-2(\xi+y_1+\dots+y_k)})_q^m = \\ = \frac{1}{(1-q^{-2})^m}\sum_{i=0}^{m}\binom{m}{i}(-q^{-2\xi})^i  \prod_{j=1}^{k}q^{-2(h-j+i)y_j}. 
    \end{eqnarray*}
    The operators $\Phi_{y_i}$ can be applied independently to each factor $q^{-2(h-j+i)y_j}$, hence applying the operator $ \Phi_{y_1}\dots\Phi_{y_k}$ to the last part of the formula, we obtain a $q$-Bernoulli polynomial. \hfill\rule{6.5pt}{6.5pt}
\end{pf}
We now use formula (\ref{Bernoulli}) to prove Lemma \ref{recurrence}.
\paragraph{Proof of formula (\ref{Polynomial->Number}).}
Note that $[\xi + y]_q = [\xi]_qq^{-2y}+[y]_q$, using this we can calculate
\begin{eqnarray*} 
\beta_{m}^{(h,k)}(\xi) &=& \Phi_{y_1}\dots\Phi_{y_k} \left( q^{-2\sum\limits_{i=1}^k(h-i)y_i}\left(\left[\xi+\sum\limits_{i=1}^ky_i\right]_q\right)^m\right) \\&=& \Phi_{y_1}\dots\Phi_{y_k} \left( q^{-2\sum\limits_{i=1}^k(h-i)y_i}\left( q^{-2\sum\limits_{i=1}^ky_i}[\xi]_q+\left[\sum\limits_{i=1}^ky_i\right]_q\right)^{\Large m} \,\right)  \\&=&\sum\limits_{j=0}^{m}([\xi]_q)^{j}\binom{m}{j}\Phi_{y_1}\dots\Phi_{y_k} \left(q^{-2\sum\limits_{i=1}^k(h-i+j)y_i}\left(\left[\sum\limits_{i=1}^ky_i\right]_q\right)^{m-j} \right)\\ &=& \sum\limits_{j=0}^{m}\binom{m}{j}([\xi]_q)^{j}\beta_{m-j}^{(h+j,k)}.
\end{eqnarray*}
    %\hfill\rule{6.5pt}{6.5pt}
\paragraph{Proof of formula (\ref{k->1}).} Note that $[x+y]_q = [x]_qq^{-2y} + [y]_q$. Moreover, a simple induction on $k$ shows that $\left[\sum\limits_{j=1}^ky_j\right]_q = \sum\limits_{j=1}^k[y_j]_q q^{-2\sum\limits_{i<j}y_i} $. Therefore,
    $$\left(\left[\sum\limits_{j=1}^ky_j\right]_q\right)^m = \left(\sum\limits_{j=1}^k[y_j]_q q^{-2\sum\limits_{i<j}y_i}\right)^m = \sum\limits_{i_1+\dots+i_k = m} \binom{m}{i_1,i_2,\dots,i_k} \prod\limits_{j=1}^{k}([y_{j}]_q)^{i_j}q^{-2\sum\limits_{i<j}i_jy_j}.$$
    To obtain the result, it remains to use linearity of the operators $\Phi_y$ and the following formula $$\Phi_{y_1}\dots\Phi_{y_k}\left(q^{-2\sum\limits_{j=1}^k(h-j)y_j}\prod\limits_{j=1}^{k}([y_{j}]_q)^{i_j}q^{-2\sum\limits_{i<j}i_jy_j} \right)= \prod\limits_{j=1}^{k}\Phi_{y_1}(([y_{j}]_q)^{i_j} q^{-2 (h-j+i_1+\dots+i_{j-1})y_i}).$$
%\hfill\rule{6.5pt}{6.5pt}

\paragraph{Proof of formula (\ref{h->1}).} 
   We have
    \begin{eqnarray*}
        \beta_{m}^{(h,1)}-(1-q^{-2})\beta_{m+1}^{(h,1)} = \Phi_y(q^{-2(h-1)y}([y]_q)^m-q^{-2(h-1)y}(1-q^{-2})([y]_q)^{m+1}) =\\= \Phi_y(q^{-2(h-1)y}([y]_q)^m (1-(1-q^{-2y}))) = \Phi_y(q^{-2hy}([y]_q)^m) = \beta_{m}^{(h+1,1)}.
    \end{eqnarray*} 
    %\hfill\rule{6.5pt}{6.5pt}

\paragraph{Proof of formulas (\ref{Clask->1}) and (\ref{ClassPol->Number}).}
    Using the obvious equality $$\sum_{m=0}^\infty B_m^{(k)}(0)\frac{t^m}{m!} = \left(\sum_{m=0}^\infty B_m^{(1)}(0)\frac{t^m}{m!}\right)^k.$$ we obtain formula (\ref{Clask->1}) for classical Bernoulli numbers. 

Formula (\ref{ClassPol->Number}) can be obtained directly from the following relation $\sum_{m=0}^\infty B_m^{(k)}(\xi)\frac{t^m}{m!}=e^{\xi t}\sum_{m=0}^\infty B_m^{(k)}(0)\frac{t^m}{m!}.$\hfill\rule{6.5pt}{6.5pt}

\end{document}